\documentclass[sigconf]{acmart}
\AtBeginDocument{%
  }

\usepackage{subcaption}

\usepackage[ruled]{algorithm2e} 

\copyrightyear{2024}
\acmYear{2024}
\setcopyright{acmlicensed}\acmConference[SIGGRAPH Conference Papers '24]{Special Interest Group on Computer Graphics and Interactive Techniques Conference Conference Papers '24}{July 27-August 1, 2024}{Denver, CO, USA}
\acmBooktitle{Special Interest Group on Computer Graphics and Interactive Techniques Conference Conference Papers '24 (SIGGRAPH Conference Papers '24), July 27-August 1, 2024, Denver, CO, USA}
\acmDOI{10.1145/3641519.3657490}
\acmISBN{979-8-4007-0525-0/24/07}
\begin{document}
\title[Preconditioned Nonlinear Conjugate Gradient Method]{Preconditioned Nonlinear Conjugate Gradient Method for Real-time Interior-point Hyperelasticity}
\author{Xing Shen}
\affiliation{%
  \institution{Fuxi AI Lab, NetEase Inc}
    \city{Hangzhou}
  \state{Zhejiang}
  \country{China}
}
\email{shenxing03@corp.netease.com}

\author{Runyuan Cai}
\affiliation{%
  \institution{Fuxi AI Lab, NetEase Inc}
    \city{Hangzhou}
  \state{Zhejiang}
  \country{China}
}
\email{cairunyuan@corp.netease.com}

\author{Mengxiao Bi}
\affiliation{%
  \institution{Fuxi AI Lab, NetEase Inc}
    \city{Hangzhou}
  \state{Zhejiang}
  \country{China}
}
\email{bimengxiao@corp.netease.com}

\author{Tangjie Lv}
\affiliation{%
  \institution{Fuxi AI Lab, NetEase Inc}
    \city{Hangzhou}
  \state{Zhejiang}
  \country{China}
}
\email{hzlvtangjie@corp.netease.com}

\begin{abstract}
The linear conjugate gradient method is widely used in physical simulation, particularly for solving large-scale linear systems derived from Newton's method. The nonlinear conjugate gradient method generalizes the conjugate gradient method to nonlinear optimization, which is extensively utilized in solving practical large-scale unconstrained optimization problems. However, it is rarely discussed in physical simulation due to the requirement of multiple vector-vector dot products. Fortunately, with the advancement of GPU-parallel acceleration techniques, it is no longer a bottleneck. In this paper, we propose a Jacobi preconditioned nonlinear conjugate gradient method for elastic deformation using interior-point methods. 
Our method is straightforward, GPU-parallelizable, and exhibits fast convergence and robustness against large time steps. The employment of the barrier function in interior-point methods necessitates continuous collision detection per iteration to obtain a penetration-free step size, which is computationally expensive and challenging to parallelize on GPUs. To address this issue, we introduce a line search strategy that deduces an appropriate step size in a single pass, eliminating the need for additional collision detection.
Furthermore, we simplify and accelerate the computations of Jacobi preconditioning and Hessian-vector product for hyperelasticity and barrier function. Our method can accurately simulate objects comprising over 100,000 tetrahedra in complex self-collision scenarios at real-time speeds. 
\end{abstract}

\begin{teaserfigure}
\includegraphics[width=0.333\textwidth]{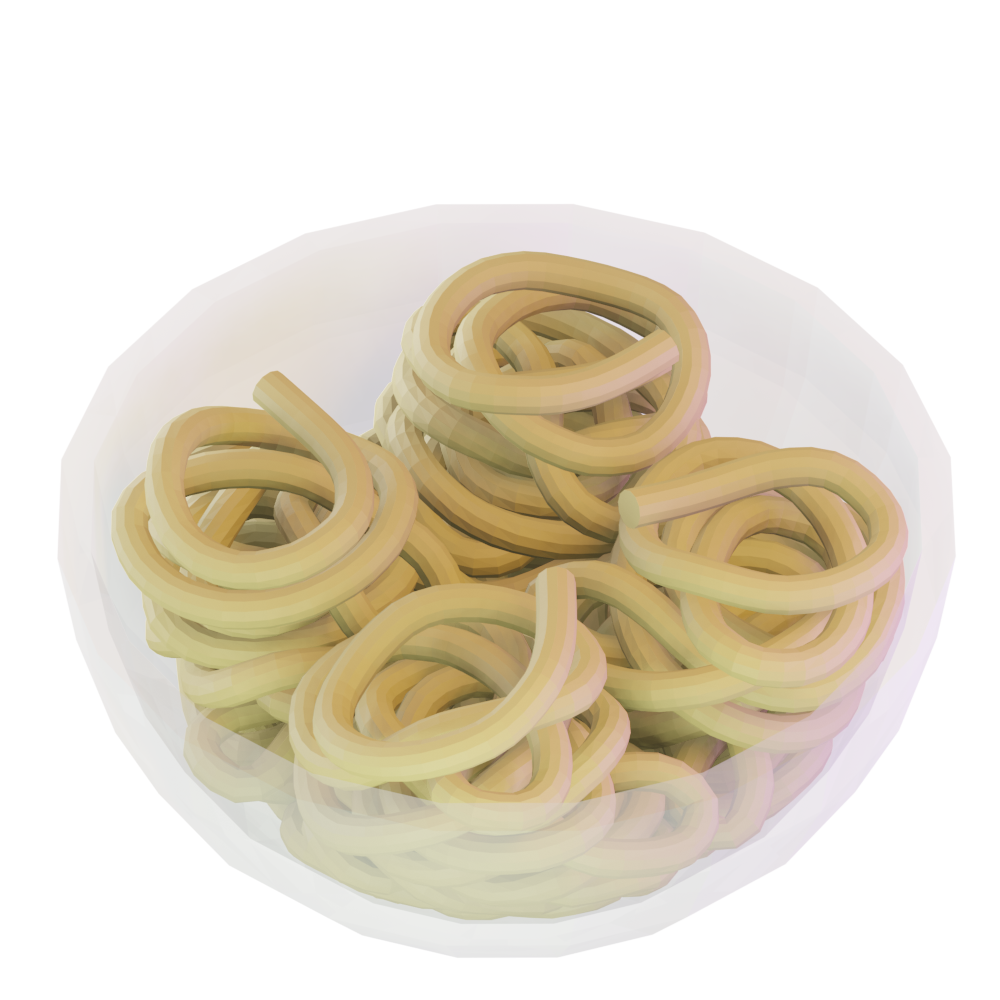} 
\includegraphics[width=0.333\textwidth]{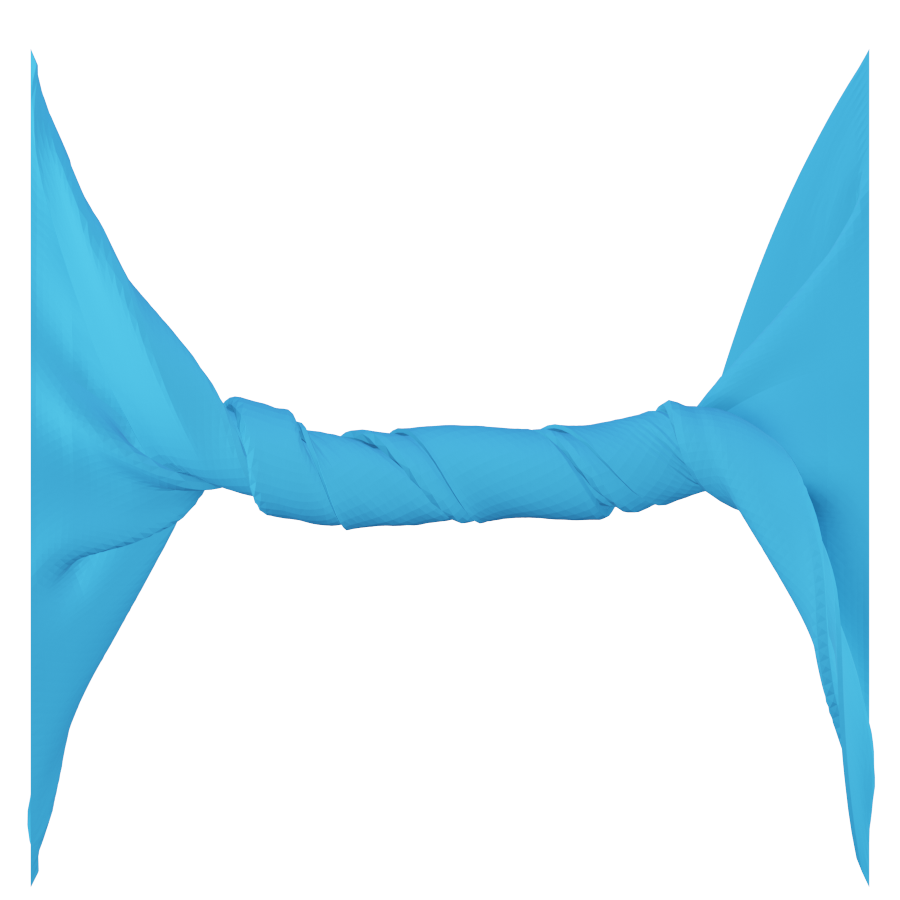}
\includegraphics[width=0.333\textwidth]{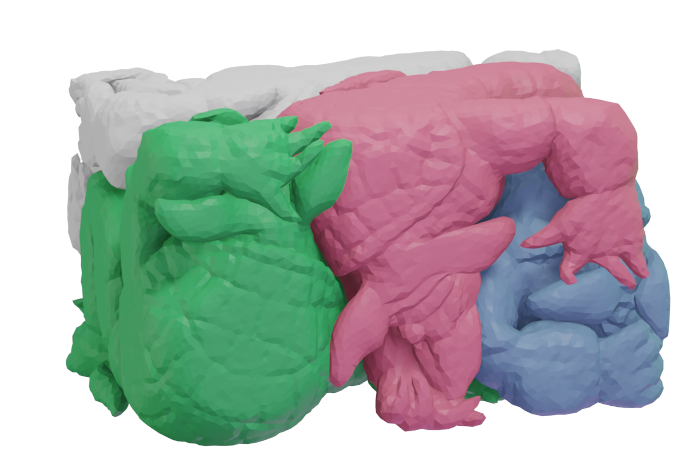}
\caption{Example simulation results involving complex self-collision scenarios.}
\end{teaserfigure}
\maketitle

\section{Introduction}
The contact simulation of elastically deformable objects is an important research topic in numerous applications, including surgical training, robotics, augmented reality/virtual reality (AR/VR), and digital fashion \cite{meier2005real,umedachi2013highly,popescu1999virtual,wang2018rule,choi2005modal}. However, accurately and robustly simulating the dynamic behavior of elastic objects with complex self-contact configurations under different external conditions is a major challenge due to the nonlinear and non-convex nature of elasticity. To address this challenge, a recent method called Incremental Potential Contact (IPC) \cite{li2020incremental} has been proposed. 

IPC incorporates the interior-point method to enable robust, accurate, and differentiable simulation of elastodynamics and contact. The method utilizes a log barrier function to approximate the inequality constraints arising from collisions and contacts, thereby converting the problem into an unconstrained optimization problem, Newton's method is then employed to solve this optimization problem. The robustness of IPC is attributed to its line search mechanism. Initially, continuous collision detection (CCD) is employed to determine the maximum step size, ensuring that objects remain penetration-free. Subsequently, a backtracking line search is implemented to achieve a decrease in the objective function.  However, the practical application of IPC is hindered by its relatively slow simulation speed. 

One main computational cost of IPC is solving the large-scale linear system derived from Newton’s method at every iteration. The linear conjugate gradient method is a common choice to solve large-scale linear systems because it can be easily parallelized on GPUs.
Since the entire optimization problem is highly nonlinear, instead of first using Newton's method for optimization and then solving the corresponding linear system with the conjugate gradient method, why not directly use the nonlinear conjugate gradient method to solve the entire optimization problem?

In paper \cite{wang2016descent}, nonlinear conjugate gradient 
is compared with the proposed Jacobi preconditioned gradient method with Chebyshev acceleration.
Although the nonlinear conjugate gradient method demonstrates better convergence than the proposed method, it exhibits poorer overall performance due to the computation cost of the vector-vector dot product. {The author reports a time consumption of 0.41ms for a vector-vector dot product over a 15K-vertices model, utilizing the NVIDIA GeForce GTX TITAN X GPU and CUDA thrust library. Fortunately, with advancements in GPU hardware and parallel techniques, the computational cost of the dot product has significantly decreased. By leveraging the Taichi language \cite{hu2019taichi}, MeshTaichi package \cite{MeshTaichi} and NVIDIA RTX 4090 GPU, the vector-vector dot product calculation now requires only 0.012 ms for the same model, and performing 4 dot products within one for-loop requires just 0.016 ms.} With the significantly reduced time consumption of vector-vector dot product, the nonlinear conjugate gradient method becomes a promising approach to discuss. 

Nonlinear conjugate gradient(NCG) method is a well-established approach for unconstrained optimization. It demonstrates a convergence speed comparable to that of the Quasi-Newton method and is particularly suitable for GPU parallel acceleration. Numerous variations of the algorithm have been proposed to improve convergence for different optimization problems \cite{andrei2020nonlinear}. However, in the field of physics simulation, researchers typically rely on the classic Fletcher-Reeves and Polak-Ribière algorithms   \cite{PR,fletcher1964function}. The essence of the nonlinear conjugate gradient method lies in its strategy for calculating the search direction, which primarily involves performing several vector-vector dot products. With our implementation, the computational cost of vector-vector dot product is relatively low, and the differences in computational time per iteration of different algorithms are negligible. Therefore, by finding a conjugate gradient algorithm that converges faster, we can directly accelerate the entire process of physical simulation. 

The line search is another significant computational burden in IPC. Firstly, the computational cost of CCD is relatively high and it is challenging to accelerate through GPU parallelization. Secondly, despite CCD providing an upper bound for the step size, a backtracking line search is still necessary for convergence. All of these operations aim to determine an appropriate step size. Therefore, if we can directly deduce a suitable step size from the available information, it would eliminate the need for collision detection and significantly speed up the algorithm. 

In this paper, we propose a preconditioned nonlinear conjugate gradient (PNCG) method for interior-point hyperelasticity and demonstrate that Dai-Kou \cite{DK} conjugate gradient algorithm has the fastest convergence rate. Our PNCG method is straightforward, GPU-parallelizable, and exhibits fast convergence and robustness against large time steps. We introduce a line search strategy to deduce an appropriate step size in one pass and achieve penetration-free simulation in complex scenarios without any additional collision detection module. Computational costs associated with the Hessian matrix are typically significant, so we simplify and parallelize the computation of the Hessian matrix for hyperelasticity and barrier function separately. The entire algorithm is parallelized on the GPU using the Taichi programming language and the MeshTaichi package. Leveraging the benefits of fast convergence optimization algorithms, improved line search strategies, simplified Hessian computations, and efficient GPU parallelization, our method enables real-time simulations of objects with over 100,000 tetrahedra in complex self-collision scenarios throughout the entire simulation. The source code is available at \url{https://github.com/Xingbaji/PNCG_IPC/}.

\section{Related work}
The simulation of elastically deformable objects has been an important graphics research topic since the 1980s \cite{terzopoulos1987elastically, terzopoulos1988constraints,terzopoulos1988deformable}. The main objective of deformable simulation is to accurately reproduce the behavior of real-world materials in a digital environment.
Modeling, theoretical analysis, and numerical simulation of elastodynamics and contact have been extensively discussed in  \cite{bergou2008discrete,kane1999finite,kaufman2008staggered,liu2017quasi,teran2005robust, verschoor2019efficient,li2020incremental}. Due to the highly nonlinear and non-convex nature of elasticity, as well as the typically non-smooth behavior of contact, achieving robust and accurate simulation of these effects remains a significant challenge. 

Recently, Li et al. \cite{li2020incremental} proposed the Incremental Potential Contact (IPC) method for robust, accurate, and differentiable elastodynamics and contact simulations. This method has the capability to consistently produce high-quality results for codimensional solids across a wide range of scenarios, ensuring interpenetration-free simulations. This method has been successfully applied in various areas, including co-dimensional simulation \cite{li2020codimensional}, rigid body simulation \cite{ferguson2021intersection,lan2022affine}, embedded FEM  \cite{zhao2022barrier}, FEM-MPM coupling \cite{li2022bfemp}, and geometric modeling \cite{fang2021guaranteed}.

Although IPC has numerous advantages, its primary limitation lies in its computational speed. The IPC method utilizes logarithmic barrier functions to approximate inequality constraints induced by collisions and contacts. It converts the overall variational optimization into an unconstrained optimization problem and utilizes the Newton's method to solve it. 
A CCD-based line search is followed to ensure all the primitives are intersection-free before any displacement update is committed. The Newton solver and the CCD-based line search are both expensive and hard to implement GPU acceleration.     
Recent efforts have focused on accelerating IPC through the use of reduced-order models \cite{lan2021medial}, projective dynamics \cite{lan2022penetration,li2023subspace}, block coordinate descent  \cite{li2023second}, and time splitting \cite{wang2023fast,xie2023contact}.
However, enabling IPC method to perform comprehensive Finite Element Analysis (FEA) at real-time speeds remains a significant challenge.

The conjugate gradient (CG) method has been a classic algorithm for solving unconstrained optimization problems since the 1950s \cite{nocedal1999numerical}.
 In the field of computer graphics, linear conjugate gradient algorithms, especially Fletcher-Reeves  \cite{fletcher1964function}, are usually employed to solve the large linear system, due to the matrix-free property and suitability for parallel execution on GPUs. Recently, several new variants of CG method  \cite{DK,HZ,DESCON,andrei2020nonlinear} are proposed and exhibit better convergence rate than Fletcher-Reeves algorithm. However, it is rarely discussed in the field of physical simulation.
The conjugate gradient algorithms are usually sensitive to the step size, and the line search strategy is the crucial point in conjugate gradient algorithms. Numerous line search strategies have been proposed for various conjugate gradient algorithms, such as \cite{HZ,DK,DESCON,liu2018several,andrei2009acceleration}. However, due to the limited number of experimental results on CG algorithms with over 100K variables, it remains a significant challenge to determine an appropriate step size.

\section{Method}
We first provide a brief overview of the formulation of elastodynamics using interior-point method \cite{mehrotra1992implementation}. By employing the implicit Euler method, the variational optimization of the elastic simulation is formulated as:
\begin{equation}
\mathbf{x}^{t+1}=\arg \min _\mathbf{x} \frac{1}{2}(\mathbf{x}-\tilde{\mathbf{x}})^{\top} \mathbf{M}(\mathbf{x}-\tilde{\mathbf{x}})+h^2 \Psi(\mathbf{x}) \quad s.t. \; h_i(\mathbf{x})\geq 0.
\end{equation}
where $\mathbf{x}^{t+1}$ denotes the positions of all vertices within a three-dimensional tetrahedral meshed model at time step $t+1$. The term $\tilde{\mathbf{x}} = \mathbf{x}^{t} + h\mathbf{v}^{t} + h^2\mathbf{M}^{-1}\mathbf{f}_{ext}$
,  where $h$ denotes time step size, $\mathbf{v}$ represents the velocity of vertices, $\mathbf{f}_{ext}$ is the external force, and $\mathbf{M}$ refers to the mass matrix. The first term $ \frac{1}{2}(\mathbf{x}-\tilde{\mathbf{x}})^{\top} \mathbf{M} (\mathbf{x}-\tilde{\mathbf{x}})$
denotes the inertia potential. The hyperelastic energy $\Psi(\mathbf{x})$ is utilized to quantify the magnitude of deformation. The inequality constraints set $C$, given by $h_i(\mathbf{x})\geq 0$, ensures that the simulation is free from both inter- and intra-model intersections.

Incremental potential contact \cite{li2020incremental} is an implementation of the interior-point method that approximates constraints $h_i(\mathbf{x})\geq 0$ with a barrier function:
$$\kappa \sum_{k \in C} b\left(d_k(\mathbf{x})\right).$$ 
The barrier function of IPC is defined as
\begin{equation}\label{barrier}
b\left(d_k(\mathbf{x})\right)= \begin{cases}-\left(d_k-\hat{d}\right)^2 \log \left(\frac{d_k}{\hat{d}}\right), & 0<d_k<\hat{d}, \\ 0, & d_k \geq \hat{d},\end{cases}
\end{equation}
where the hyperparameter $\hat{d}$ controls the threshold of collision repulsion, and $\kappa$ determines the magnitude of collision repulsion.

Consequently, IPC converts the original problem into the following unconstrained optimization problem:
\begin{align}
\mathbf{x}^{t+1}&=\arg \min _\mathbf{x} E(\mathbf{x}),\\ E(\mathbf{x}) &= \frac{1}{2}(\mathbf{x}-\tilde{\mathbf{x}})^{\top} \mathbf{M}(\mathbf{x}-\tilde{\mathbf{x}})+h^2 \Psi(\mathbf{x}) + \kappa \sum_{k \in C} b\left(d_k(\mathbf{x})\right)\label{E}.
\end{align}
The Hessian matrix  of $E(\mathbf{x})$ is given by:
\begin{equation}\label{Hessian}
\mathbf{H} = \mathbf{M} + h^2\frac{\partial^2 \Psi}{\partial \mathbf{x}^2} + \kappa \sum_{k \in C} \frac{\partial^2 b}{\partial \mathbf{x}^2}.
\end{equation}
In this paper, we need to compute the diagonal matrix of $\mathbf{H}$ as the preconditioner and $\mathbf{p}^{\top}\mathbf{H}\mathbf{p}$ in the line search, where $\mathbf{p}$ represents the search direction in NCG algorithms. Note that the computational cost of these two terms will significantly impact the speed of algorithm. Therefore, we aim to optimize the computation of the Hessian matrix for $\Psi$ and $b$ separately. 

\subsection{Hyperelasticity}
In this paper, we employ following invariants described in Chapter 5 of \cite{kim2022dynamic} to formulate the energy density function for a hyperelastic material.
$$I_1=\operatorname{tr}(\mathbf{S}), \quad I_2=\operatorname{tr}\left(\mathbf{F}^T \mathbf{F}\right), \quad I_3=\operatorname{det}( \mathbf{F}).$$
where $\mathbf{F}$ denotes the deformation gradient and $\mathbf{S}$ is derived from Polar decomposition $\mathbf{F} = \mathbf{R}\mathbf{S} $.

The elastic behavior of a deformable body can be characterized in terms
of a hyperelastic energy density $\Psi(I_1,I_2,I_3)$, such as the Neo-Hookean elasticity \cite{bonet_wood_2008}:
$$\Psi_{\mathrm{NH}}=\frac{\mu}{2}\left(I_2-3\right)-\mu \log \left(I_3\right)+\frac{\lambda}{2}\left(\log \left(I_3\right)\right)^2.$$

The gradients and Hessians of each invariant  with respect to $\mathbf{F}$ can be obtained as follows:
$$\begin{aligned}\mathbf{g}_1 &=\operatorname{vec}(\mathbf{R}),&\mathbf{H}_1=&\sum_{i=0}^2 \lambda_i  \mathbf{q}_i\mathbf{q}_i^{\top},\\ \mathbf{g}_2 &= 2 \operatorname{vec}(\mathbf{F}), & \mathbf{H}_2= & \;2\mathbf{I}_{9 \times 9}, \\\mathbf{g}_3 &= \operatorname{vec}\left(\left[\mathbf{f}_1 \times \mathbf{f}_2, \mathbf{f}_2 \times \mathbf{f}_0, \mathbf{f}_0 \times \mathbf{f}_1 \right]\right), & \mathbf{H}_3= &\left[\begin{array}{rrr} \mathbf{0} & -\hat{\mathbf{f}}_2 & \hat{\mathbf{f}}_1 \\\hat{\mathbf{f}}_2 & \mathbf{0} & -\hat{\mathbf{f}}_0 \\-\hat{\mathbf{f}}_1 & \hat{\mathbf{f}}_0 & \mathbf{0}\end{array}\right] ,
\end{aligned}$$
where $\mathbf{f}_0$, $\mathbf{f}_1$ and $\mathbf{f}_2$ represent the columns of $\mathbf{F}$. $\mathbf{I}$ represents the identity matrix, and $\mathbf{0}$ represents the zero matrix. Vectorization $vec(\cdot)$ denotes column-wise flattening of a matrix into a vector, the symbol $\hat{\mathbf{f}}_i$ indicates the cross-product matrix out of the vector $\mathbf{f}_i$.
The detailed definitions of vector $\mathbf{q}_i$, and matrices $\mathbf{H}_1$ and $\mathbf{H}_3$ are provided in the supplemental document.

Now we can derive the Hessian matrix of hyperelasticity:
$$\begin{aligned} \frac{\partial^2 \Psi}{\partial \mathbf{x}^2}=& \operatorname{vec}\left(\frac{\partial \mathbf{F}}{\partial \mathbf{x}}\right)^{\top}  \operatorname{vec}\left(\frac{\partial^2 \Psi}{\partial \mathbf{F}^2}\right) \operatorname{vec}\left(\frac{\partial \mathbf{F}}{\partial \mathbf{x}}\right) \\ 
=& \operatorname{vec}\left(\frac{\partial \mathbf{F}}{\partial \mathbf{x}}\right)^{\top}\left( \sum_{i=1}^3( \frac{\partial^2 \Psi}{\partial I_i^2} \mathbf{g}_i \mathbf{g}_i^{\top}+\frac{\partial \Psi}{\partial I_i} \mathbf{H}_i)\right)\operatorname{vec}\left(\frac{\partial \mathbf{F}}{\partial \mathbf{x}}\right),
\end{aligned}$$
where $\operatorname{vec}\left(\frac{\partial \mathbf{F}}{\partial \mathbf{x}}\right) \in \mathcal{R}^{9\times 12}$ and $\operatorname{vec}\left(\frac{\partial^2 \Psi}{\partial \mathbf{F}^2}\right) \in \mathcal{R}^{9\times 9}$. 

{
Hence we can split the Hessian matrix into six terms,
$$\frac{\partial^2 \Psi}{\partial \mathbf{x}^2}= \sum_{i=1}^3\left( \frac{\partial^2 \Psi}{\partial I_i^2} \mathbf{h}_i + \frac{\partial \Psi}{\partial I_{i}} \mathbf{h}_{i+3}\right),$$
where
$$
\begin{aligned}\mathbf{h}_i=&\operatorname{vec}\left(\frac{\partial \mathbf{F}}{\partial \mathbf{x}}\right)^{\top} \mathbf{g}_i \mathbf{g}_i^{\top} \operatorname{vec}\left(\frac{\partial \mathbf{F}}{\partial \mathbf{x}}\right), \; \text{for}\; i=1,2,3,\\
\mathbf{h}_{i} &=\operatorname{vec}\left(\frac{\partial \mathbf{F}}{\partial \mathbf{x}}\right)^{\top}\mathbf{H}_{i} \operatorname{vec}\left(\frac{\partial \mathbf{F}}{\partial \mathbf{x}}\right), \; \text{for}\; i = 4,5,6. \end{aligned}$$
So we have
$$\begin{aligned}\overline{diag}\left(\frac{\partial^2 \Psi}{\partial \mathbf{x}^2}\right) =& \sum_{i=1}^3\left( \frac{\partial^2 \Psi}{\partial I_i^2} \overline{diag}(\mathbf{h}_i) + \frac{\partial \Psi}{\partial I_{i}} \overline{diag}(\mathbf{h}_{i+3})\right),\\
\mathbf{p}^{\top}\frac{\partial^2 \Psi}{\partial \mathbf{x}^2}\mathbf{p} =& \sum_{i=1}^3\left( \frac{\partial^2 \Psi}{\partial I_i^2} \mathbf{p}^{\top}\mathbf{h}_i\mathbf{p} + \frac{\partial \Psi}{\partial I_{i}} \mathbf{p}^{\top}\mathbf{h}_{i+3}\mathbf{p}\right),
\end{aligned}$$
where the notation $\overline{diag}(A)$ to represent the vector composed of the diagonal elements of matrix $A$, and vector $\mathbf{p} \in \mathcal{R}^{12\times 1}$ is the search direction of NCG algorithm. }

In this approach, the computation of $\frac{\partial^2 \Psi}{\partial \mathbf{x}^2}$ is divided into six parts, each of which can be efficiently calculated with minimal floating-point operations (FLOPs). {In the supplemental document, we provide detailed instructions on computing $\overline{diag}(\mathbf{h_i})$ and $\mathbf{p}^{\top}\mathbf{h}_{i}\mathbf{p}$.}
For the Neo-Hookean elasticity model, the calculation of $\overline{diag}\left(\frac{\partial^2 \Psi}{\partial \mathbf{x}^2}\right)$ and $\mathbf{p}^{\top}\frac{\partial^2 \Psi}{\partial \mathbf{x}^2}\mathbf{p}$ requires only 143 and 314 FLOPs, respectively. However, if we first compute $\frac{\partial^2 \Psi}{\partial \mathbf{F}^2}$ and then compute $\frac{\partial^2 \Psi}{\partial \mathbf{x}^2}$ using two matrix-matrix products, it necessitates 624 FLOPs and 3590 FLOPs, respectively.
Since these two terms will be computed in parallel for every tetrahedral element, the reduced computational cost with our method would be significant. Additionally, the computation of gradient and $\overline{diag}(\mathbf{H})$ can be included in the same for-loop, further enhancing the computational speed.

Compared to the implementation in \cite{li2023second}, which computes Hessian-vector product with accelerated complex-step finite difference (CSFD) \cite{luo2019accelerated}, our method is more concise in form and faster in speed. Moreover, our method can handle any hyperelastic material and does not rely on simplifications of the material models such as projective dynamics  \cite{PD,lan2022penetration} and position-based dynamics \cite{macklin2016xpbd}.

\subsection{Incremental Potential Contact}
Inspired by Chapter 14 of \cite{kim2022dynamic}, we make slight modifications to the definition of the distance in IPC between primitives in order to improve parallel acceleration.

The distance between a point $\mathbf{x}_0$ and a triangle $T = (\mathbf{x}_1,\mathbf{x}_2,\mathbf{x}_3)$ is defined as follows:
\begin{align}\label{PT}
d^{PT} = & \| \mathbf{x}_0 - \mathbf{x}_b \| = \min_{\alpha,\beta,\gamma}\| \mathbf{x}_0 - (\alpha \mathbf{x}_1 + \beta \mathbf{x}_2 + \gamma \mathbf{x}_3) \|, \\  & s.t.\; \alpha,\beta, \gamma \geq 0, \; \alpha + \beta + \gamma = 1,    \nonumber
\end{align}
where  $\mathbf{x}_b = \alpha \mathbf{x}_1 + \beta \mathbf{x}_2 + \gamma \mathbf{x}_3$ represents the point within the triangle that is closest to the point $\mathbf{x}_0$. 
Similarly, the distance between edges $\mathbf{x}_0\text{-}\mathbf{x}_1$ and $\mathbf{x}_2\text{-}\mathbf{x}_3$ can be defined as follows:
\begin{equation}\label{EE}
\begin{aligned} 
{d}^{EE} = & \| \mathbf{x}_a - \mathbf{x}_b \| = \min_{a_0,a_1,b_0,b_1}\| (a_0 \mathbf{x}_0 + a_1 \mathbf{x}_1)  - (b_0 \mathbf{x}_2 + b_1 \mathbf{x}_3) \|, \\  & s.t.\; a_0,a_1,b_0,b_1 \geq 0, \; a_0 +a_1 = 1,\; b_0 + b_1 = 1,
\end{aligned}
\end{equation}
where $\mathbf{x}_a= a_0 \mathbf{x}_0 + a_1 \mathbf{x}_1$ and $\mathbf{x}_b = b_0 \mathbf{x}_2 + b_1 \mathbf{x}_3$  are the closest points between the edges. By combining Equation (\ref{PT}) and Equation(\ref{EE}), we can reformulate the distance computation as follows:
{
\begin{align}
{d} = & \| \mathbf{t} \| = \min_{c_0,c_1,c_2,c_3}\| c_0 \mathbf{x}_0 + c_1 \mathbf{x}_1 + c_2 \mathbf{x}_2 + c_3 \mathbf{x}_3 \|, \\ 
  & s.t. \; c_0 = 1,\; c_1,c_2,c_3\leq 0,\; c_1+c_2+c_3=-1 \text{ for PT primitives},   \nonumber\\ 
  & c_0,c_1\geq 0,\; c_2,c_3\leq 0,\; c_0+c_1=1,\; c_2+c_3 = -1 \text{ for EE primitives}.
  \nonumber
\end{align}
}
Distance $d$ is represented as the norm of  $\mathbf{t} \in \mathcal{R}^{3\times 1}$, which is the linear combination of $\mathbf{x}_0,\mathbf{x}_1,\mathbf{x}_2,\mathbf{x}_3$. If either $c_1,c_2$ or $c_3$ is equal to 0 in PT primitives,  the point-edge(PE) distance is derived. If exactly two of $c_1,c_2$ or $c_3$ are equal to $0$, the point-point(PP) distance is derived. Similarly, the PP and PE distances can be derived from EE primitives. 

Hence, the gradient and Hessian of $t$ are 
$$\frac{\partial \mathbf{t}}{\partial \mathbf{x}}=\left[\begin{array}{llll} c_0 \mathbf{I}_{3 \times 3} & c_1 \mathbf{I}_{3 \times 3} & c_2 \mathbf{I}_{3 \times 3} & c_3\mathbf{I}_{3 \times 3}\end{array}\right]^{\top} \in \mathcal{R}^{12 \times 3}, \space \frac{\partial^2 \mathbf{t}}{\partial \mathbf{x}^2}=\mathbf{0},$$
It is worth noting that the definition of distance remains the same as the original formulation presented in \cite{li2020incremental}, the only difference is that we introduce $\mathbf{t}$ to simplify the computation.

Therefore, the gradient and Hessian of barrier function are
\begin{align}
&\frac{\partial b }{\partial \mathbf{x}} = \frac{\partial b }{\partial d} \frac{\partial d }{\partial \mathbf{t}} \frac{\partial \mathbf{t} }{\partial \mathbf{x}} = \frac{\partial b }{\partial d} \frac{1}{d} \frac{\partial \mathbf{t} }{\partial \mathbf{x}} \mathbf{t}, \label{eq8}   \\ 
&\frac{\partial^2 b }{\partial \mathbf{x}^2}  = (\frac{\partial^2 b }{\partial d^2} \frac{1}{d^2} - \frac{\partial b }{\partial d} \frac{1}{d^3}) (\frac{\partial \mathbf{t} }{\partial \mathbf{x}} \mathbf{t})(\frac{\partial \mathbf{t} }{\partial \mathbf{x}} \mathbf{t})^{\top}  + \frac{\partial b }{\partial d} \frac{1}{d} (\frac{\partial \mathbf{t} }{\partial \mathbf{x}} )(\frac{\partial \mathbf{t} }{\partial \mathbf{x}} )^{\top}. \label{eq9}
\end{align} 
Because $(\frac{\partial^2 b }{\partial d^2} \frac{1}{d^2} - \frac{\partial b }{\partial d} \frac{1}{d^3})$ and $\frac{\partial b }{\partial d} \frac{1}{d}$ are scalar, we focus on the matrix terms.  Therefore, we obtain 
\begin{align}
&\overline{diag} ((\frac{\partial \mathbf{t} }{\partial \mathbf{x}}\ )(\frac{\partial \mathbf{t} }{\partial \mathbf{x}} )^{\top}) = [c_0,c_0,c_0,c_1,c_1,c_1,c_2,c_2,c_2,c_3,c_3,c_3]^{\circ 2}, \label{eq10}\\
&\overline{diag} ( (\frac{\partial \mathbf{t} }{\partial \mathbf{x}} \mathbf{t} )(\frac{\partial \mathbf{t} }{\partial \mathbf{x}} \mathbf{t} )^{\top}) = [c_0 \mathbf{t}^{\top}, c_1 \mathbf{t}^{\top}, c_2 \mathbf{t}^{\top}, c_3 \mathbf{t}^{\top}]^{\circ 2}, \label{eq11}
\end{align}
where the Hadamard power $\circ 2$ represents element-wise square of the vector.
Thus, we can efficiently compute the diagonal matrix without the need to calculate the entire Hessian matrix.

Given a vector $\mathbf{p} = [\mathbf{p}_0,\mathbf{p}_1,\mathbf{p}_2,\mathbf{p}_3]^{\top} \in \mathcal{R}^{12\times 1}$ , where $\mathbf{p}_i \in \mathcal{R}^{1\times 3}$, we have 
\begin{align}
& \mathbf{p}^{\top}(\frac{\partial \mathbf{t} }{\partial \mathbf{x}} t)(\frac{\partial \mathbf{t} }{\partial \mathbf{x}} t)^{\top}\mathbf{p} = (\mathbf{p}^{\top}\frac{\partial \mathbf{t} }{\partial \mathbf{x}} t)^2, \label{eq12}\\
& \mathbf{p}^{\top} (\frac{\partial \mathbf{t} }{\partial \mathbf{x}} )(\frac{\partial t}{\partial \mathbf{x}} )^{\top} \mathbf{p} = \|\mathbf{p}^{\top}\frac{\partial \mathbf{t} }{\partial \mathbf{x}}  \|_2^2, \label{eq13}
\end{align}
where the computation of $\mathbf{p}^{\top}\frac{\partial t}{\partial \mathbf{x}} = \sum_{i=0}^3 c_i\mathbf{p}_i$ is straightforward. This approach significantly reduces the computational cost compared to initially calculating a $12 \times 12$ matrix and subsequently performing two matrix multiplications.

With equation (\ref{eq9}),(\ref{eq10}) and (\ref{eq11}), we can simplify the computation of $\overline{diag}(\frac{\partial^2 b}{\partial \mathbf{x}^2})$. Similarly, with equation (\ref{eq9}), (\ref{eq12}) and (\ref{eq13}), we can simplify the computation of $\mathbf{p}^{\top} \frac{\partial^2 b}{\partial \mathbf{x}^2} \mathbf{p}$. Additionally, the procedure for identifying the constraints derived from PT and EE primitives closely follows that of  \cite{li2020incremental}, with the only difference being the computation and storage of $\mathbf{t}$ and the corresponding coefficients $c_0, c_1, c_2, c_3$. The impact of these differences on computational time and memory usage is negligible. Furthermore, by merging the definition and corresponding computation of PT distance and EE distance, we can handle all the constraints within a single for-loop. Thus the entire process can be parallelized for acceleration.

Through the above implementation, the computation of $\overline{diag}(\mathbf{H})$ and $\mathbf{p}^{\top} \mathbf{H} \mathbf{p}$ only requires two for-loops, one for all the tetrahedral elements and another for all the contact constraints. Both of these loops are totally parallelizable and very fast.
\section{Nonlinear Conjugate Gradient Algorithm}
As illustrated in Algorithm \ref{alg}, the pipeline of the proposed preconditioned nonlinear conjugate gradient method is illustrated. In this section, we describe the techniques employed in our method.

\subsection{Conjugate Gradient Algorithm} 
Through a series of experiments, we find that the Dai-Kou (DK) algorithm \cite{DK} exhibits superior performance in terms of robustness and convergence speed. The formula of DK conjugate gradient algorithm is 
\begin{align}
\mathbf{p}_{k+1}&=- \mathbf{g}_{k+1}+{\beta}_k^{D K} \mathbf{p}_k, \\
{\beta}_k^{D K}&=\frac{\mathbf{g}_{k+1}^T \mathbf{y}_k}{\mathbf{y}_k^T \mathbf{p}_k}-\frac{\mathbf{y}_k^T  \mathbf{y}_k}{\mathbf{y}_k^T \mathbf{p}_k} \frac{\mathbf{p}_k^T \mathbf{g}_{k+1}}{\mathbf{y}_k^T \mathbf{p}_k},\end{align}
where $\mathbf{p}_k$ is the conjugate gradient direction at iteration $k$,
 $\mathbf{g}_k$ is the gradient of objective function, and  $\mathbf{y}_k = \mathbf{g}_{k+1} - \mathbf{g}_k$. 

The formula for $\beta_{DK}$ is derived by determining the conjugate gradient direction that is closest to the direction of the scaled memoryless BFGS method (SSML-BFGS), and the numerical experiments in \cite{DK} showed that DK algorithm performs more efficiently and robustly than the SSML-BFGS method. The proof of convergence and theoretical analysis in \cite{DK} also illustrate its efficiency.

\subsection{Preconditioning}  
Preconditioning is a widely employed technique to accelerate the conjugate gradient algorithms by reducing the condition number of the problem. As discussed in Section 3, the diagonal elements of Hessian of IPC and hyperelastic can be easily computed and accelerated with GPU parallel technique. Therefore, we adopt Jacobi preconditioning matrix $\mathbf{P} = diag(\mathbf{H})^{-1}$, which is the inverse of the diagonal matrix of the Hessian matrix. So we can derive the preconditioned DK algorithm:
\begin{align}
\mathbf{p}_{k+1}&=-\mathbf{P}_{k+1} \mathbf{g}_{k+1}+{\beta}_k^{D K} \mathbf{p}_k \\
{\beta}_k^{D K}&=\frac{\mathbf{g}_{k+1}^T \mathbf{P}_{k+1} \mathbf{y}_k}{\mathbf{y}_k^T \mathbf{p}_k}-\frac{\mathbf{y}_k^T \mathbf{P}_{k+1} \mathbf{y}_k}{\mathbf{y}_k^T \mathbf{p}_k} \frac{\mathbf{p}_k^T \mathbf{g}_{k+1}}{\mathbf{y}_k^T \mathbf{p}_k},
\end{align}
where the computation of $\beta^{DK}$ requires four vector-vector dot products.

For vertices with collision repulsion, as the distance between pairs of primitives approaches $0$, the value of $b$ exhibits a significant increase due to the logarithmic term of Equation (\ref{barrier}). This increase leads to larger values of $diag(\mathbf{H})$, resulting in a smaller scale of search direction $p$. Conversely, vertices without collision repulsion exhibit a larger scale of search direction, leading to a similar convergence rate among all vertices. As a result, the application of preconditioning can effectively alleviate the "locking" issue mentioned in \cite{li2023second}.

In Section 3, we have optimized the computational complexity of the diagonal Hessian matrix to the same level as computing the gradient. Additionally, the gradient and preconditioner can be evaluated in the same for-loop. Consequently, despite evaluating $\mathbf{P}$ at each iteration, the computational cost is still acceptable, as opposed to evaluating the Hessian matrix once every $M$ ($M>>1$) iterations in \cite{wang2016descent}, which may lead to artifacts or penetration in complex contact scenarios and exhibit slightly poorer convergence.

\subsection{Line Search} 
Conjugate gradient methods are usually sensitive to the step size. A step size that is too small would decrease the convergence speed, while a step size that is too large may result in penetration or artifacts. Therefore, it is crucial to compute an appropriate step size.
Nevertheless, elasticity is highly nonlinear and non-convex, and the  Hessian matrix of IPC may change drastically between two iterations due to the logarithmic barrier function. As a result, neither the interpolation method nor the Barzilai-Borwein method \cite{barzilai1988two} can determine an appropriate step size. {Therefore, we adopt Newton's method for line search, which secures an optimal step size by leveraging a quadratic approximation of the objective function. By utilizing the Taylor series, we obtain
\begin{equation}
\label{Taylor}
E(\mathbf{x}+\alpha \mathbf{p}) \approx E(\mathbf{x}) + \alpha \mathbf{g}^{\top}\mathbf{p} + \frac{\alpha^2}{2} \mathbf{p}^{\top} \mathbf{H} \mathbf{p} 
\end{equation}
Hence, 
\begin{equation}
\frac{\partial E(\mathbf{x}+\alpha \mathbf{p})}{\partial \alpha} \approx  \mathbf{g}^{\top}\mathbf{p} + \alpha \mathbf{p}^{\top} \mathbf{H} \mathbf{p} 
\end{equation}
The function $E(\mathbf{x}+\alpha \mathbf{p})$ is approximately minimized by setting the gradient with respect to $\alpha$ to be zero, so we have}
\begin{equation}
\label{linesearch}
\widebar{\alpha} = - \frac{\mathbf{g}_{k+1}^{\top}\mathbf{p}_{k+1}}{\mathbf{p}_{k+1}^{\top}\mathbf{H}_{k+1}\mathbf{p}_{k+1}}. 
\end{equation}
In Section 3, we have significantly accelerated the computation of $\mathbf{p}_{k+1}^{\top}\mathbf{H}_{k+1}\mathbf{p}_{k+1}$, and the computation of the vector-vector dot product $\mathbf{g}_{k+1}^{\top}\mathbf{p}_{k+1}$ is inexpensive. Hence the overall computational cost is affordable. Moreover, Newton's method exhibits quadratic convergence, allowing for the direct determination of an appropriate step size.

Recall that the barrier function treats the distances between primitives larger than $\hat{d}$ as $0$. In complex self-collision scenarios, when the displacement of vertices exceeds $\hat{d}$ during iterations, primitives may penetrate before contact constraints are detected. {Inspired by the "CFL-Inspired Culling of CCD" technique proposed in \cite{li2020incremental}, we impose a restriction on the step size to ensure that the maximum displacement of vertices is less than $\frac{\hat{d}}{2}$. The upper bound of step size is 
\begin{equation}\label{upper}
\alpha_{\text{upper}} = \frac{\hat{d}}{2 \| \mathbf{p}_{k+1}\|_{\infty}},
\end{equation}
so the final step size is given as:
\begin{equation}\label{alpha}
\alpha = \min (\alpha_{\text{upper}}, \widebar{\alpha})
\end{equation}
Let us denote the optimal step size as $\alpha^{*}$. When $\alpha^{*} < \alpha_{\text{upper}}$, the original IPC method initiates with  $\alpha_{\text{upper}}$ and performs backtracking line search to determine a proper step size. In contrast, our method directly obtains an appropriate step size $\widebar{\alpha}$. When $\alpha^{*} > \alpha_{\text{upper}}$, IPC employs CCD to establish a higher upper bound. Since CCD computation is computationally expensive, and each iteration of the PNCG method is computationally efficient, we prefer to utilize the conservative step size $\alpha_{\text{upper}}$ and perform additional iterations. Moreover, as the PNCG method operates as a first-order optimization technique, its search direction is comparatively less precise than that of Newton's method, typically resulting in a smaller optimal step size  $\alpha^{*}$ for each iteration. Consequently, the scenario where $\widebar{\alpha} > \alpha_{\text{upper}}$ is infrequent, primarily occurring during the initial few iterations. Limiting the step size in the early stages helps prevent optimization divergence.}

\begin{algorithm}
\caption{Preconditioned Nonlinear Conjugate Gradient}\label{alg}
$\mathbf{x}_0 \gets x^{t}$\\
$\tilde{x} = x^{t} + hv^{t} + h^2M^{-1}f_{ext}$\\
\For{$k = 0$ to IterMax }{
 $C\gets \text{ComputeConstraintSet}(x,\hat{d})$\\
 $\mathbf{g}_{k+1}, \mathbf{P}_{k+1}\gets \text{ComputeGradientAndPreconditioning}(x,\tilde{x},C)$\\
 \uIf{$k=0$}
 {
  ${\beta}_k \gets 0$
  }
 \Else{
 ${\beta}_k \gets \text{ComputeBeta}(\mathbf{g}_{k+1},\mathbf{g}_k,\mathbf{p}_{k+1},\mathbf{p}_k)$\\
 }
 $\mathbf{p}_{k+1} \gets -\mathbf{P}_{k+1} \mathbf{g}_{k+1} + \beta_k \mathbf{p}_k $\\
 $\alpha \gets \min ( \frac{\hat{d}}{2 \| \mathbf{p}_{k+1}\|_{\infty}}, - \frac{\mathbf{g}_{k+1}^{\top}\mathbf{p}_{k+1}}{\mathbf{p}_{k+1}^{\top}\mathbf{H}_{k+1}\mathbf{p}_{k+1} }  ).$\\
 $\mathbf{x}_{k+1} \gets \mathbf{x}_k + \alpha \mathbf{p}_{k+1} $\\
 $\Delta E \gets - \alpha \mathbf{g}_{k+1}^{\top}\mathbf{p}_{k+1} -\frac{\alpha^2}{2} \mathbf{p}_{k+1}^{\top} \mathbf{H}_{k+1} \mathbf{p}_{k+1}$\\
\If{$i=0$}{
 $\Delta E_0 \gets \Delta E$ 
 }
\If{$\Delta E < \epsilon \Delta E_0$}
{
 \textbf{break}
}
}
\end{algorithm}

\subsection{Termination Condition} 
From Equation (\ref{Taylor}), we can derive{
$$\Delta E = E(\mathbf{x}_k) - E(\mathbf{x}_{k+1}) \approx - \alpha \mathbf{g}_{k+1}^{\top}\mathbf{p}_{k+1} -\frac{\alpha^2}{2} \mathbf{p}_{k+1}^{\top} \mathbf{H}_{k+1} \mathbf{p}_{k+1}.$$}
When $\Delta E$ is small, the algorithm can hardly decrease the objective function, indicating the termination of iterations. 
Since the initial energy scale can vary significantly across different materials and situations, we employ the condition $\Delta E < \epsilon \Delta E_0$ as the termination condition. This termination condition does not require additional computation, and experiments show that it is robust for PNCG algorithm.

\begin{table*}[htbp]
\caption{Performance results. The "8 E Falling" scenario utilizes a time step size of $\Delta t=0.01$, all other examples utilize $\Delta t=0.04$. "\# tets" refers to the number of tetrahedra. {The abbreviations ARAP, FCR, NH, and SNH represent the As-Rigid-As-Possible, Fixed Corotated, Neo-Hookean, Stable Neo-Hookean\cite{smith2018stable} elasticity models, respectively.}}
\begin{tabular}{|l|l|l|l|l|l|l|l|}
\hline
Scene               & Material & \# nodes, \# tets \# faces & \# contact avg(max) & IterMax & Tolerance $\epsilon$ & avg Iters & FPS avg(min) \\ \hline
{Drag armadillo} (Fig. \ref{pics:draga}) & SNH  & 13K, 55K, 22K             & 0K (0K) & 150 & $5 \times 10^{-5}$          & 108.1         & 40.1 (25.0)   \\ \hline
{8 "E" falling (Fig. \ref{pics:E2})}           & NH     & 8K, 27K, 13K               & 2K (8K) & 50 &  $3 \times 10^{-4}$ & 32.8         & 46.6 (24.3)   \\ \hline
48 "E" falling (Fig. \ref{pics:E})           & ARAP     & 50K, 164K, 78K               & 63K (99K) & 25 & $1 \times 10^{-3}$           & 25.0         & 32.2 (27.9)   \\ \hline
Four long noodles (Fig. \ref{pics:noodles})      & ARAP  & 39K, 101K, 73K                & 13K (29K) & 25  & $1 \times 10^{-3}$            & 24.9         & 27.6 (25.8)   \\ \hline
Squeeze four armadillo (Fig. \ref{pics:armadillo}) & ARAP   & 53K, 168K, 91K             & 71K (322K) & 25 & $1 \times 10^{-3}$          & 24.1         & 28.9 (25.6)   \\ \hline
{Twist mat (Fig. \ref{pics:twist_mat})}           & FCR      & 45K, 133K, 90K               & 87K (222K) & 150 &  $1 \times 10^{-3}$          & 138         & 7.4 (5.6)   \\ \hline
{Twist four rods (Fig. \ref{pics:twist_rods})}       & FCR      & 53K, 202K, 80K               & 39K (92K) & 150  & $1 \times 10^{-3}$           & 107.5         & 10.2 (6.9)  \\ \hline
\end{tabular}
\label{table}
\end{table*}

\section{IMPLEMENTATION DETAILS}
\subsection{Larger $\mathbf{\hat{d}}$ Strategy}
In \cite{li2023second}, the authors adopt a warm start strategy by using larger $\hat{d}$, which leads to a significant reduction in the total number of iterations. Our experiments also demonstrate that a larger $\hat{d}$ enhances the convergence of the PNCG algorithm. It enables the algorithm to detect collisions earlier and utilize a smaller step size, effectively preventing penetration. Moreover, a larger $\hat{d}$ relaxes the restriction on the maximum step size in Equation (\ref{alpha}). Instead of using different customized $\hat{d}$ values for each collision pair, a constant $\hat{d}$ would be sufficient for our method. However, a larger $\hat{d}$ can cause additional self-collision repulsion at the rest pose. Therefore we propose a preprocessing stage. By setting the threshold as $1.5\hat{d}$ in the rest pose, we can first filter out the self-collision constraints, and store them with a spatial hashing table. By referring to the hashing table, we can exclude these self-collision constraints during the simulation.

Although a larger value of $\hat{d}$ would increase the number of contact constraints, the resulting increase in computational cost would be justifiable when considering the improvement in convergence speed, owing to the simplified and parallelized relevant calculations. Experimental results indicate that setting $\hat{d}$ within the range of 30\% to 70\% of the average boundary edge length is appropriate. Since collision detection operations are not included in the line search, it is essential to choose appropriate hyperparameters for the barrier function. Increasing the value of $\hat{d}$ can enhance algorithm's robustness against variations in the value of $\kappa$.

\subsection{Splitting Computation}  In scenarios involving multiple objects, such as Figure \ref{pics:E}, certain objects may undergo collisions and compression, while others remain collision-free. It is not appropriate for all objects to be assigned the same step size in such cases. Hence, it is necessary to utilize separate values of step size $\alpha$ and search direction $p$ for each object. This strategy can be easily achieved through parallel computation, minimizing the impact on computational cost. Similarly, in the case of the long noodles model discussed in Figure \ref{pics:noodles}, collisions may occur between specific points of the same object, while leaving other parts unaffected. To address this, the objects can be divided into collision and non-collision parts at the start of iterations. Then, separate values of $\alpha$ and $p$ can be assigned to each part, hence achieving higher convergence speed. More implementation details are elaborated in the supplementary material.

\section{Experiments}
Our algorithm is implemented on a desktop workstation with an NVIDIA RTX 4090 GPU and an Intel Core i9-13900X CPU, using the Taichi language for GPU acceleration with single-precision float. Tetrahedral mesh computation is accelerated with MeshTaichi \cite{MeshTaichi}, and a grid-based spatial hashing technique \cite{Muller2023matthias} is employed for broad-phase collision culling.

\subsection{Ablation study} 
\begin{figure*}[htbp]
\centering
\begin{subfigure}[t]{0.32\textwidth}
\includegraphics[width=\textwidth]{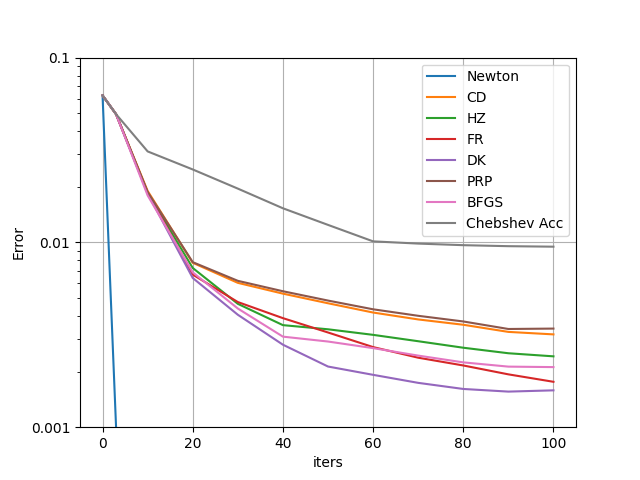}
\caption{Convergence of different methods. DK algorithm demonstrates superior performance compared to other PNCG algorithms (CD,HZ,FR,PRP), and significantly outperforms the Chebyshev accelerated preconditioned descent method. Error is evaluated using the infinity norm with respect to the ground truth.}
\label{pics:convergence}
\end{subfigure}
\hspace{\fill}
\begin{subfigure}[t]{.32\textwidth}
\includegraphics[width=\linewidth]{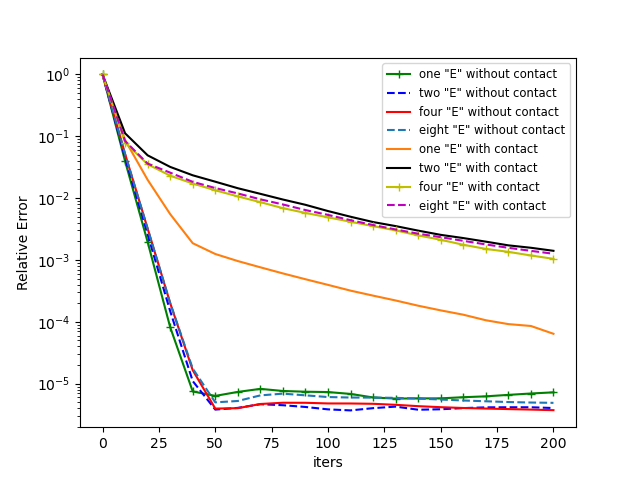}
\caption{Convergence with and without contact. In scenarios devoid of contact, an increase in the number of objects "E" has a minimal impact on convergence. In the presence of contact, the convergence speed decreases; yet, the adoption of a splitting computation strategy yields a convergence speed that is nearly identical for two, four, and eight "E".}
\label{pics:contact}
\end{subfigure}
\hspace{\fill}
\begin{subfigure}[t]{0.32\textwidth}
\includegraphics[width=\textwidth]{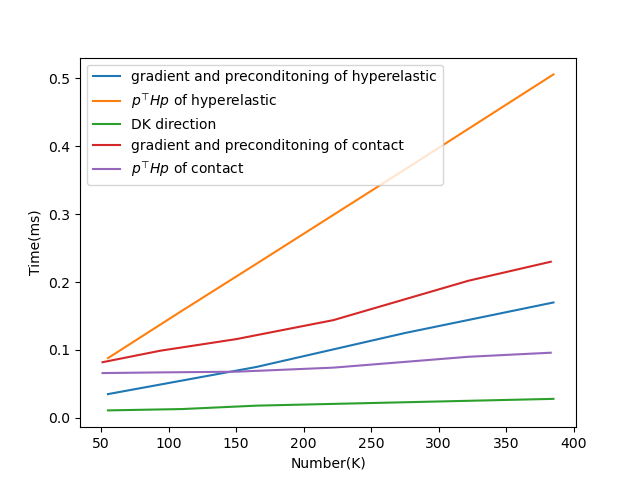}
\caption{Computational cost of hyperelastic and contact with different number of elements and constraints, measured in milliseconds. All the operation of our method shows significantly low computational costs even with a large number of elements and constraints.} 
\label{pics:time}
\end{subfigure}
\caption{Ablation study.}
\end{figure*}

In Figure \ref{pics:convergence}, we present the convergence analysis of various methods in Example \ref{pics:E}. Jacobi preconditioning is utilized for all these methods.
The ground truth is obtained by Newton-PCG with sufficient iterations. It is evident that the PNCG algorithms: CD \cite{CD}, HZ \cite{HZ}, FR \cite{fletcher1964function}, PRP \cite{polyak1969conjugate}, and DK \cite{DK}, outperform the Chebyshev-accelerated method \cite{wang2016descent}, with the DK algorithm exhibiting the highest performance among them. {The time consumption of these PNCG algorithms is  close, requiring only 1.2ms per iteration. Meanwhile, the Chebyshev-accelerated method takes slightly longer, at 1.4ms per iteration. } Furthermore, the DK method demonstrates superior performance compared to the SSML-BFGS method with Jacobi preconditioning.

Additionally, although Newton's method exhibits the best convergence, it requires more than 400ms for each iteration, significantly larger than the PNCG method. Given that the height of "E" is 1.1, simulated results are considered sufficiently realistic when the error is less than 0.01. The line search strategy for Newton's method adheres to \cite{li2020incremental}, while the other methods employ our line search strategy for consistency. The detailed formulas for these methods are provided in the supplemental document.

Figure \ref{pics:contact} compares the convergence rates of our method with and without contact scenarios. The results demonstrate that our method exhibits rapid convergence without contact. In scenarios involving contact, the convergence rate becomes slower. Because of the splitting computation strategy, the convergence rates for simulations involving two, four, and eight  "E" are nearly identical.

\begin{figure*}[htbp]
\centering
\begin{subfigure}{0.32\textwidth}
\centering
\includegraphics[width=\textwidth]{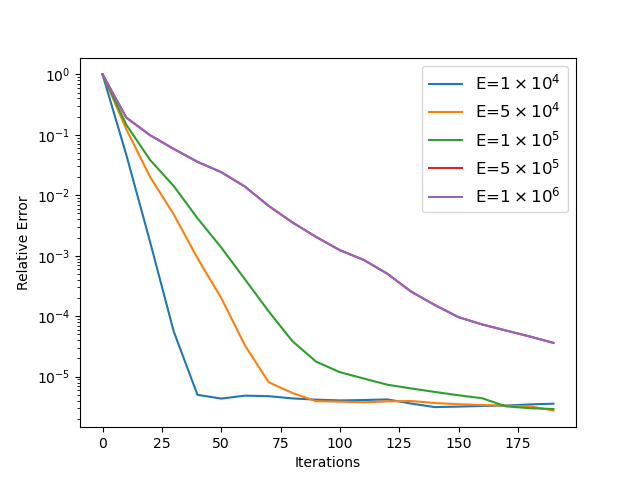}
\caption{The Impact of Young’s Modulus on Convergence.}
\label{pics:stiff}
\end{subfigure}
\hspace{\fill}
\begin{subfigure}{0.32\textwidth}
\centering
\includegraphics[width=\textwidth]{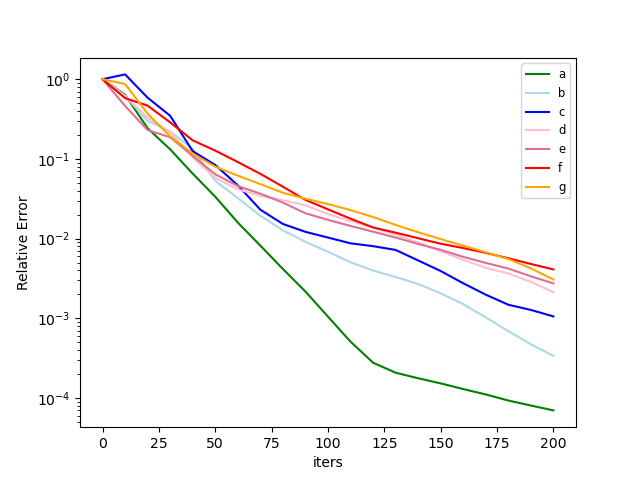}
\caption{The convergence rates for different magnitudes of deformation.}
\label{pics:drag}
\end{subfigure}
\hspace{\fill}
\begin{subfigure}{0.32\textwidth}
\centering
\includegraphics[width=\textwidth]{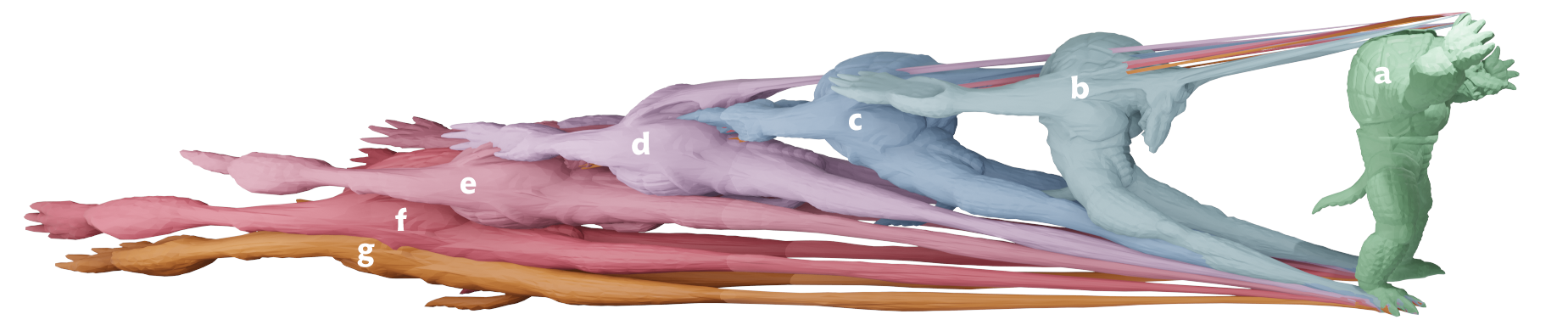}
\caption{Drag armadillo.
The markings and colors of the armadillo correspond one-to-one with Figure (b).}
\label{pics:draga}
\end{subfigure}
\caption{Additional ablation study.}
\end{figure*}

Figure \ref{pics:draga} tests our method's ability to handle extreme deformations by dragging the armadillo with large external forces in a collision-free scenario. Figure \ref{pics:drag} compares corresponding convergence rates for different deformation magnitudes. Despite substantial external forces and extreme deformations, our method exhibits satisfactory convergence, and requirs only 0.25 ms per iteration.

{Figure \ref{pics:stiff} demonstrates the free fall of one "E" object under gravity $9.8 m/s^2$ using Neo-Hookean hyperelasticity and $\Delta t =0.01$, investigating the influence of stiffness on convergence. As Young's modulus increases from 1e4 to 1e6, the convergence rate decreases, indicating our method's unsuitability for simulating stiff objects.}

\subsection{Computational Cost} Figure \ref{pics:time} shows the computational costs of hyperelasticity and barrier contact, varying with the number of tetrahedral elements and contact constraints. The simplified computations outlined in Section 3 significantly reduce these costs. Notably, the DK search direction cost is extremely low at approximately 0.01 ms. For a 220K tetrahedra model on an NVIDIA GeForce RTX 2060 GPU, the DK search direction cost is still low at 0.021 ms.

\subsection{Examples} 
\begin{figure*}[htp]
\centering
\begin{subfigure}{0.24\textwidth}
    \includegraphics[width=\textwidth]{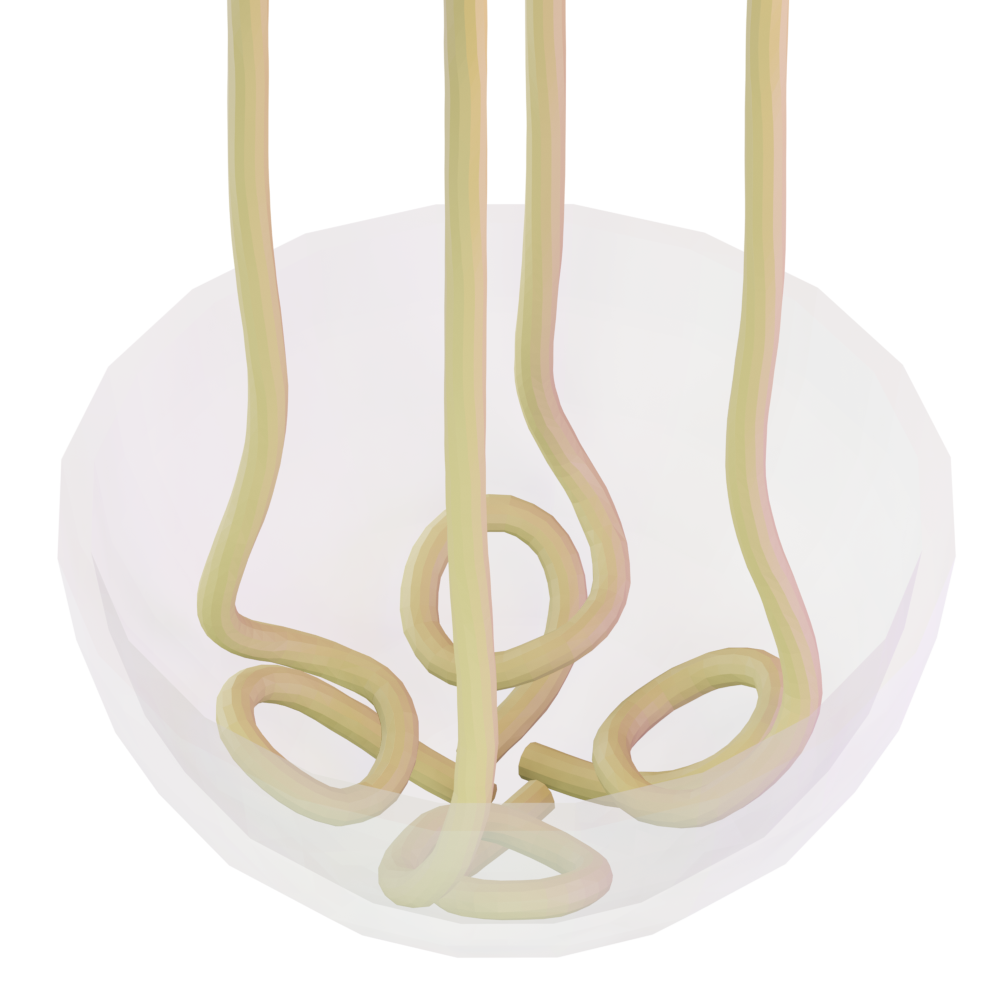}
\end{subfigure}
\begin{subfigure}{0.24\textwidth}
    \includegraphics[width=\textwidth]{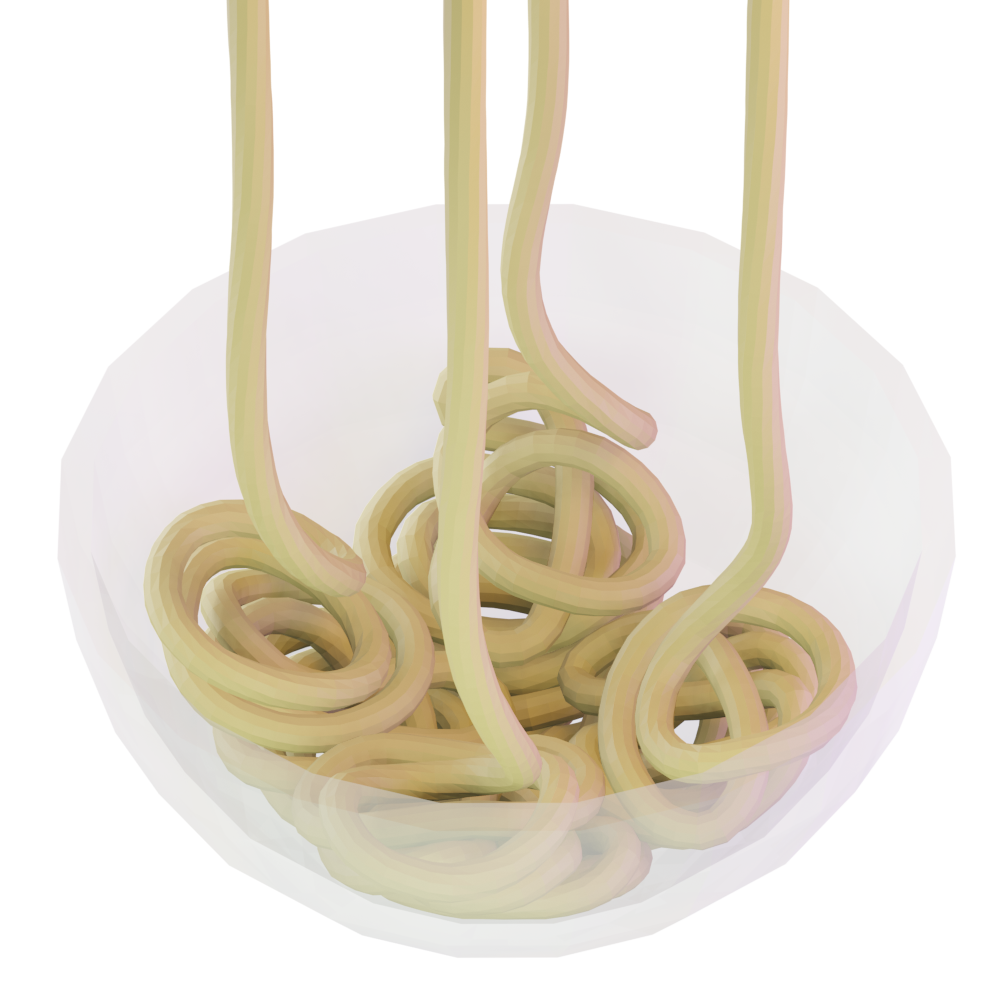}
\end{subfigure}
\begin{subfigure}{0.24\textwidth}
    \includegraphics[width=\textwidth]{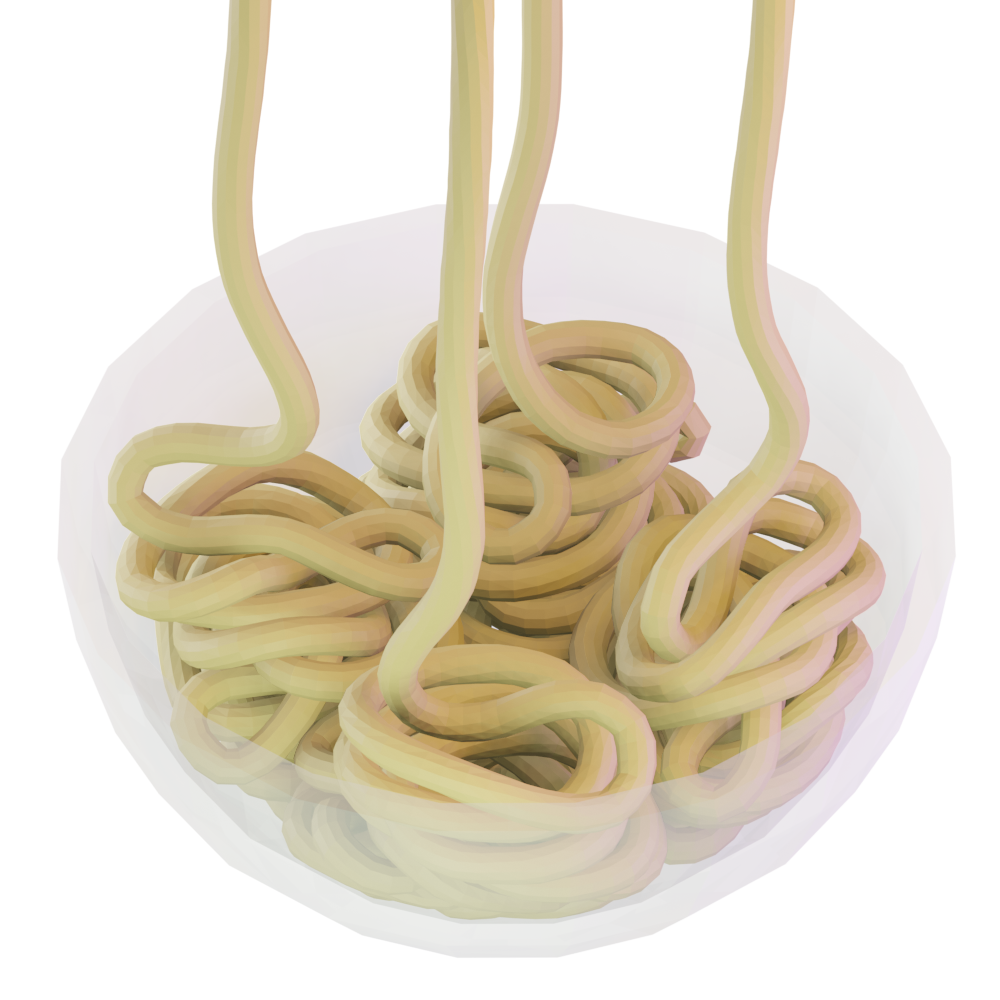}
\end{subfigure}
\begin{subfigure}{0.24\textwidth}
    \includegraphics[width=\textwidth]{imgs/noodles/1530_2.png}
\end{subfigure}
\caption{Long noodles in a bowl.}
\label{pics:noodles}
\end{figure*}

\begin{figure*}[htp]
\centering
\begin{subfigure}{0.245\textwidth}
    \includegraphics[width=\textwidth]{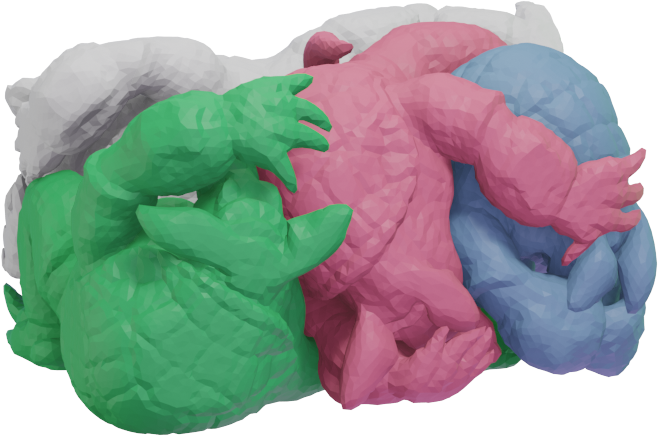}
\end{subfigure}
\begin{subfigure}{0.245\textwidth}
    \includegraphics[width=\textwidth]{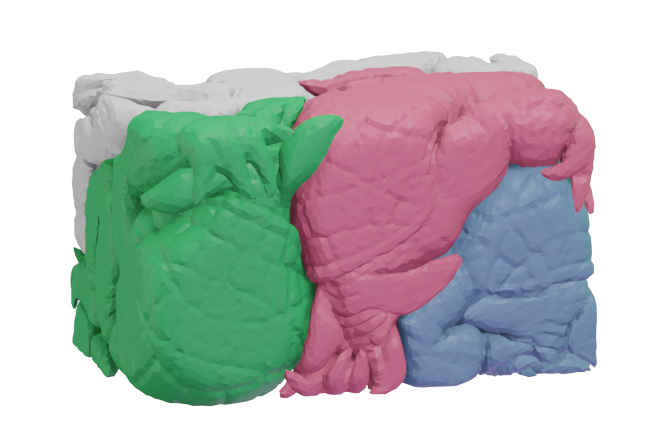}
\end{subfigure}
\begin{subfigure}{0.245\textwidth}
    \includegraphics[width=\textwidth]{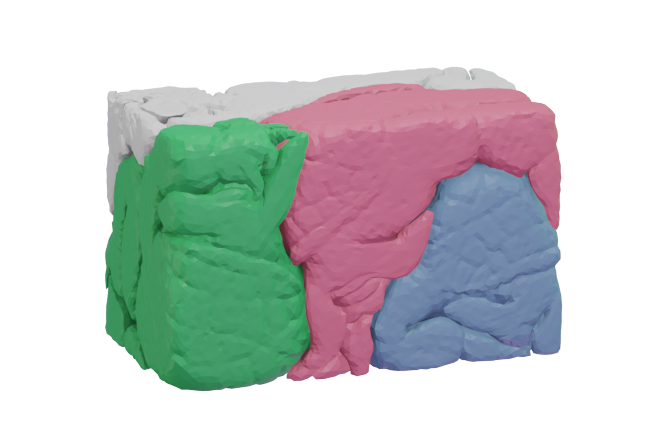}
\end{subfigure}
\begin{subfigure}{0.245\textwidth}
    \includegraphics[width=\textwidth]{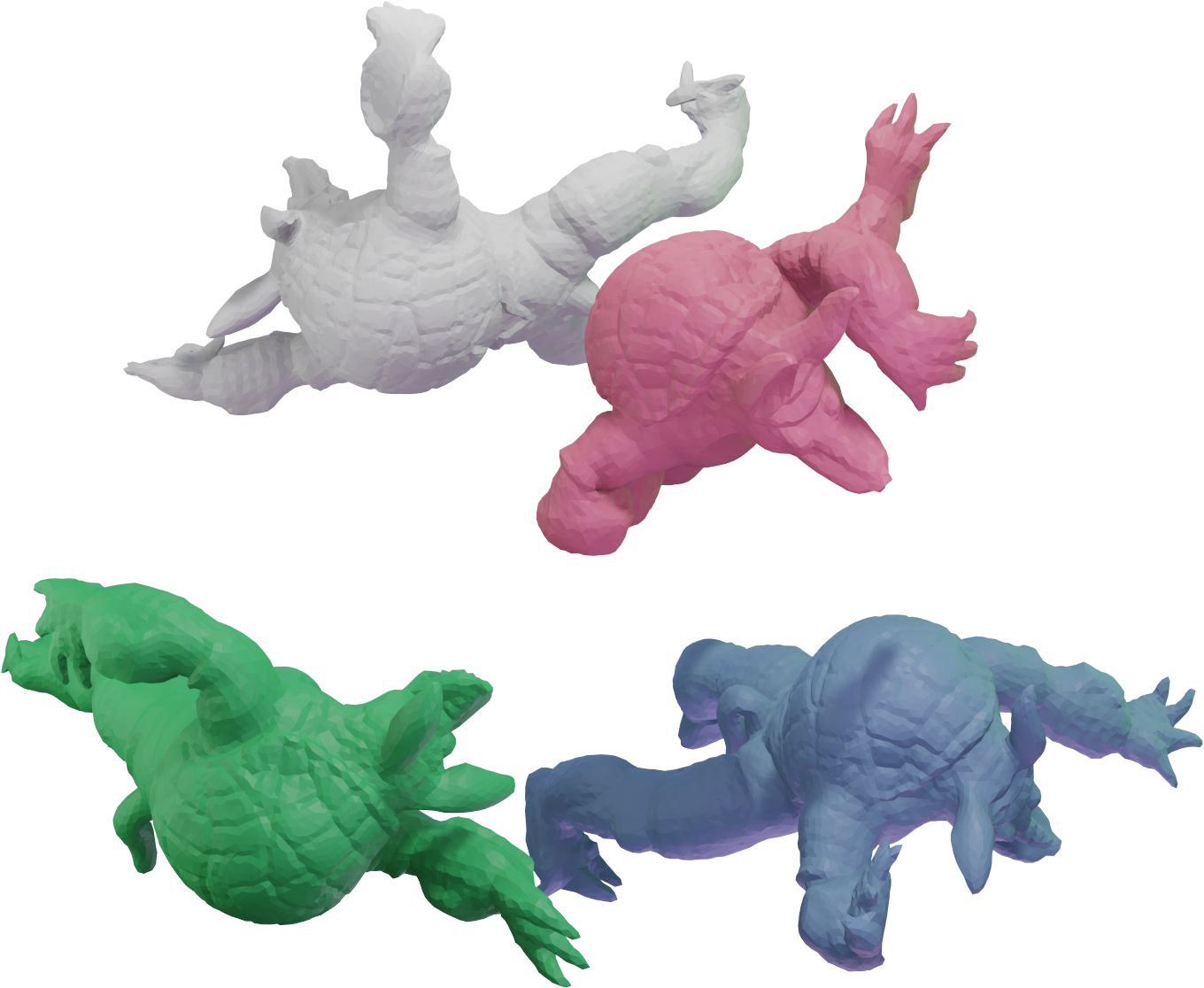}
\end{subfigure}
\caption{Squeeze and release of four armadillo models. We simultaneously compress the four small models from all six sides and then release them.}
\label{pics:armadillo}
\end{figure*}

\begin{figure*}[htp]
\centering
\begin{subfigure}{.245\textwidth}
  \centering
  \includegraphics[width=\linewidth]{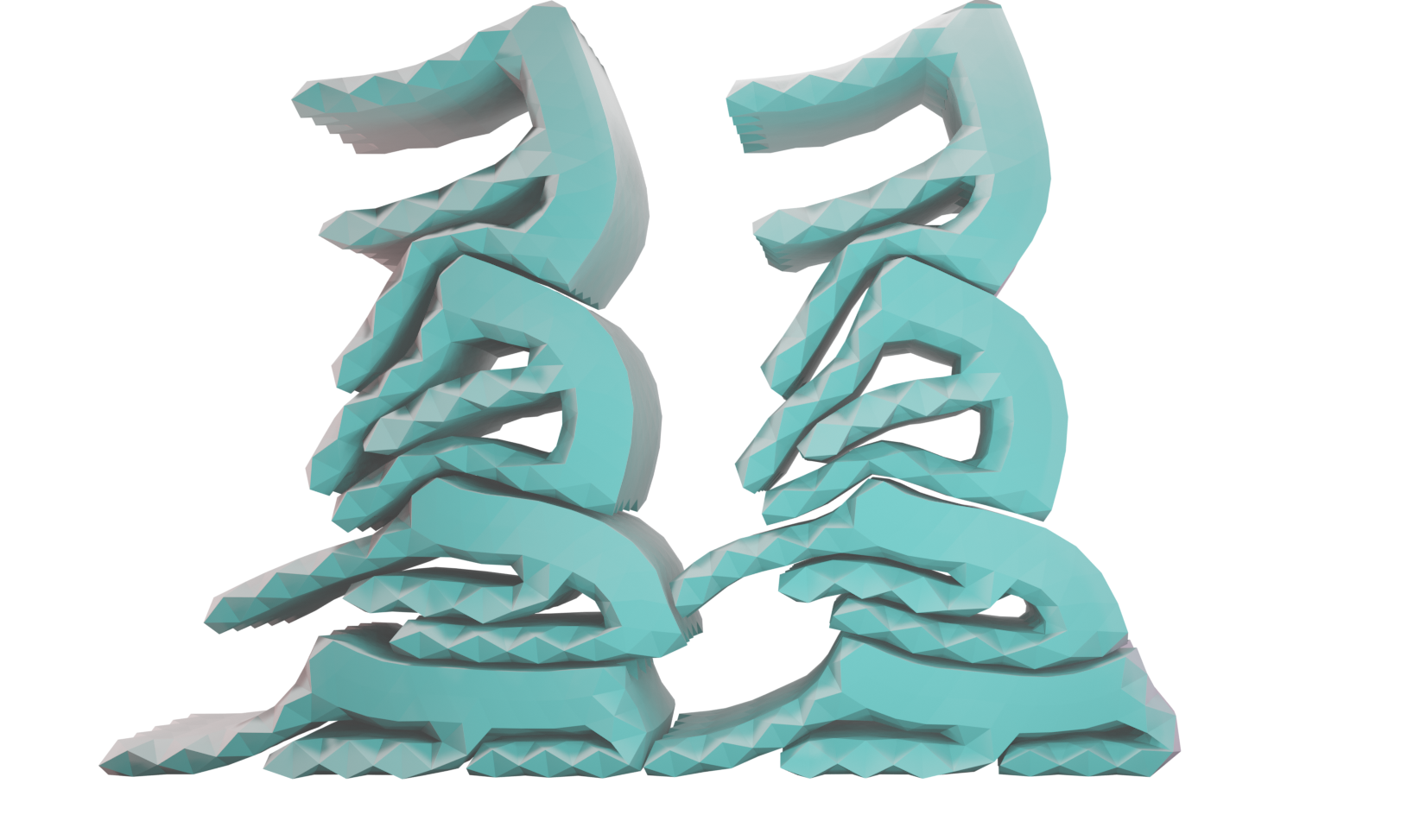}
\end{subfigure}%
\begin{subfigure}{.245\textwidth}
  \centering
  \includegraphics[width=\linewidth]{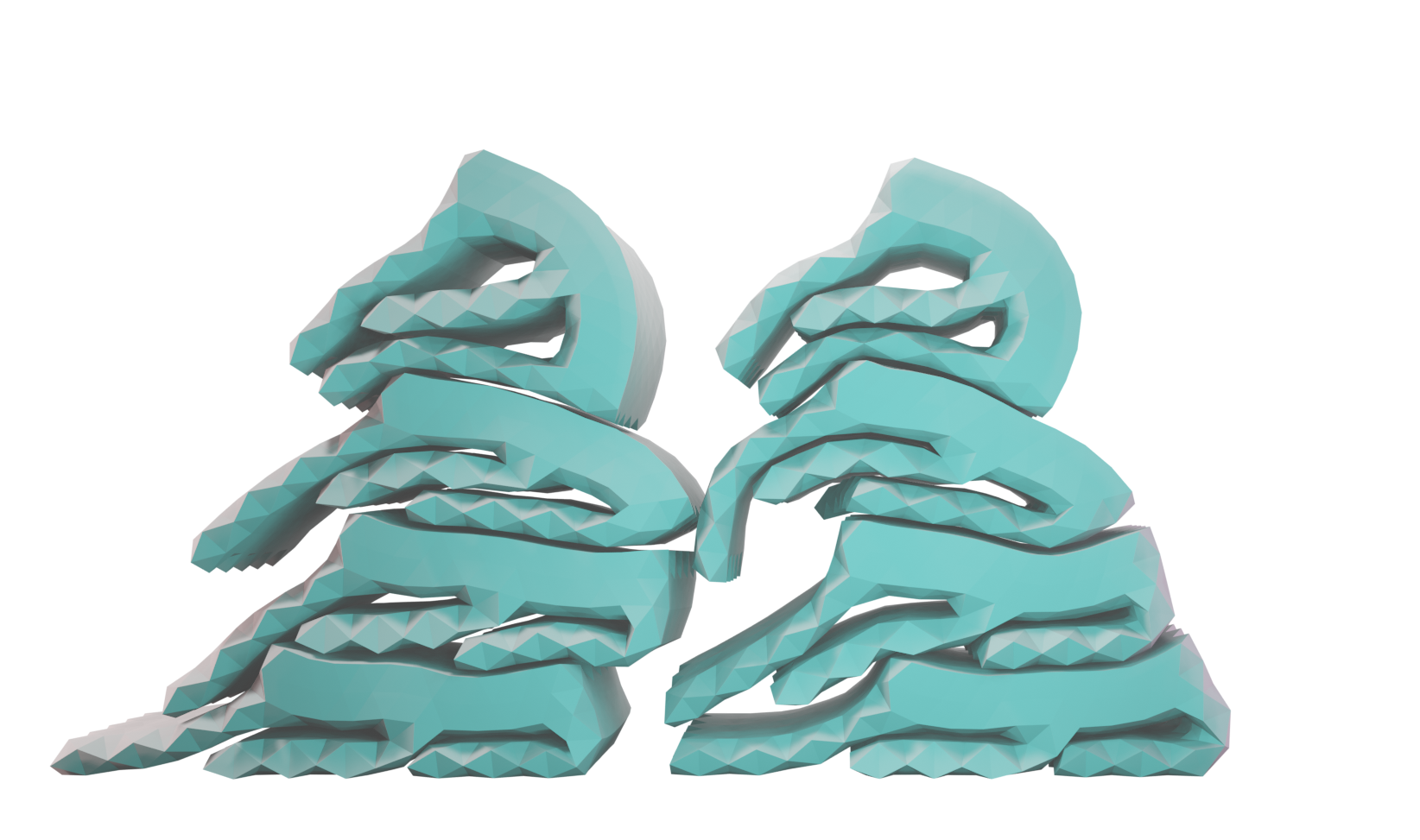}
\end{subfigure}
\begin{subfigure}{.245\textwidth}
  \centering
  \includegraphics[width=\linewidth]{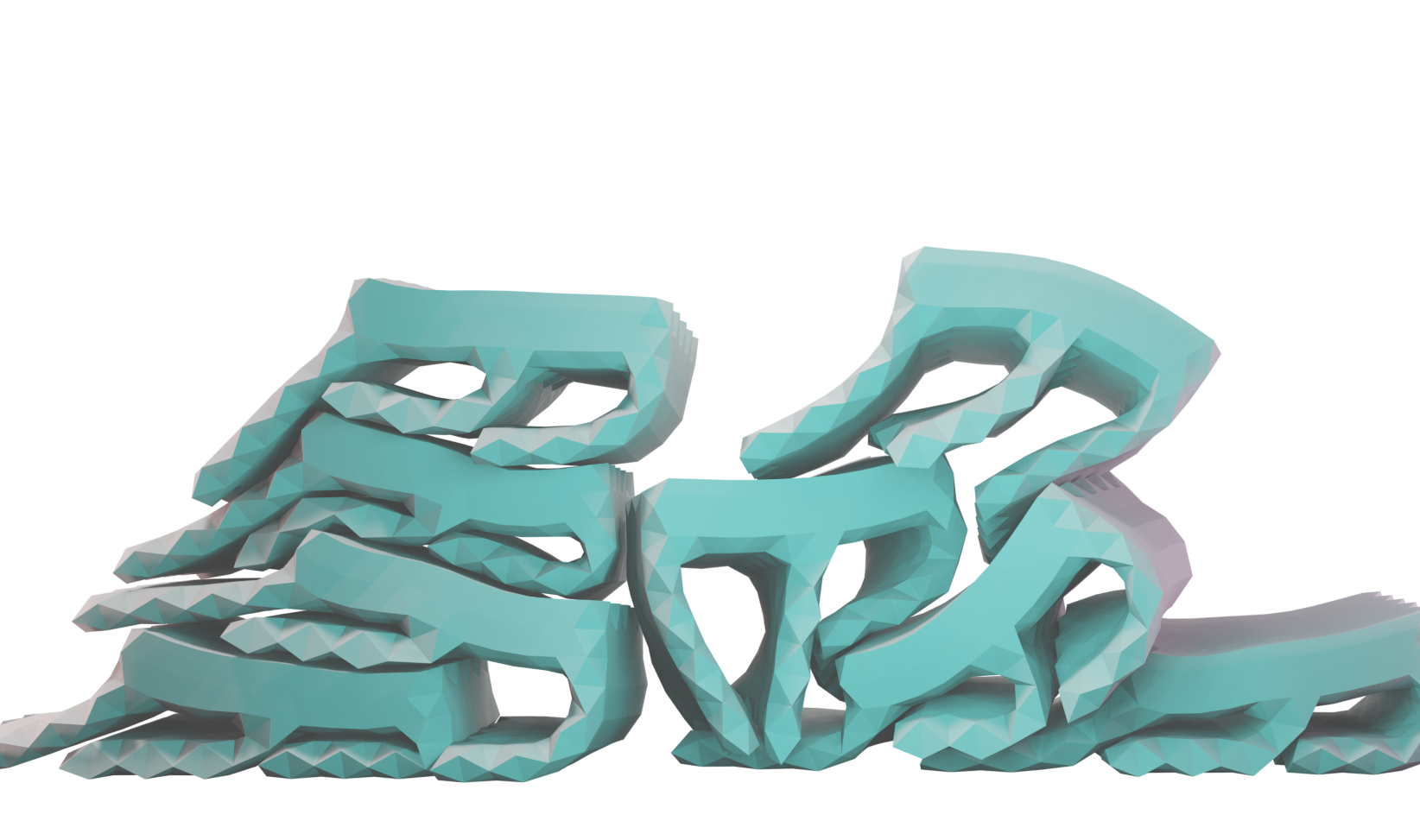}
\end{subfigure}
\begin{subfigure}{.245\textwidth}
  \centering
  \includegraphics[width=\linewidth]{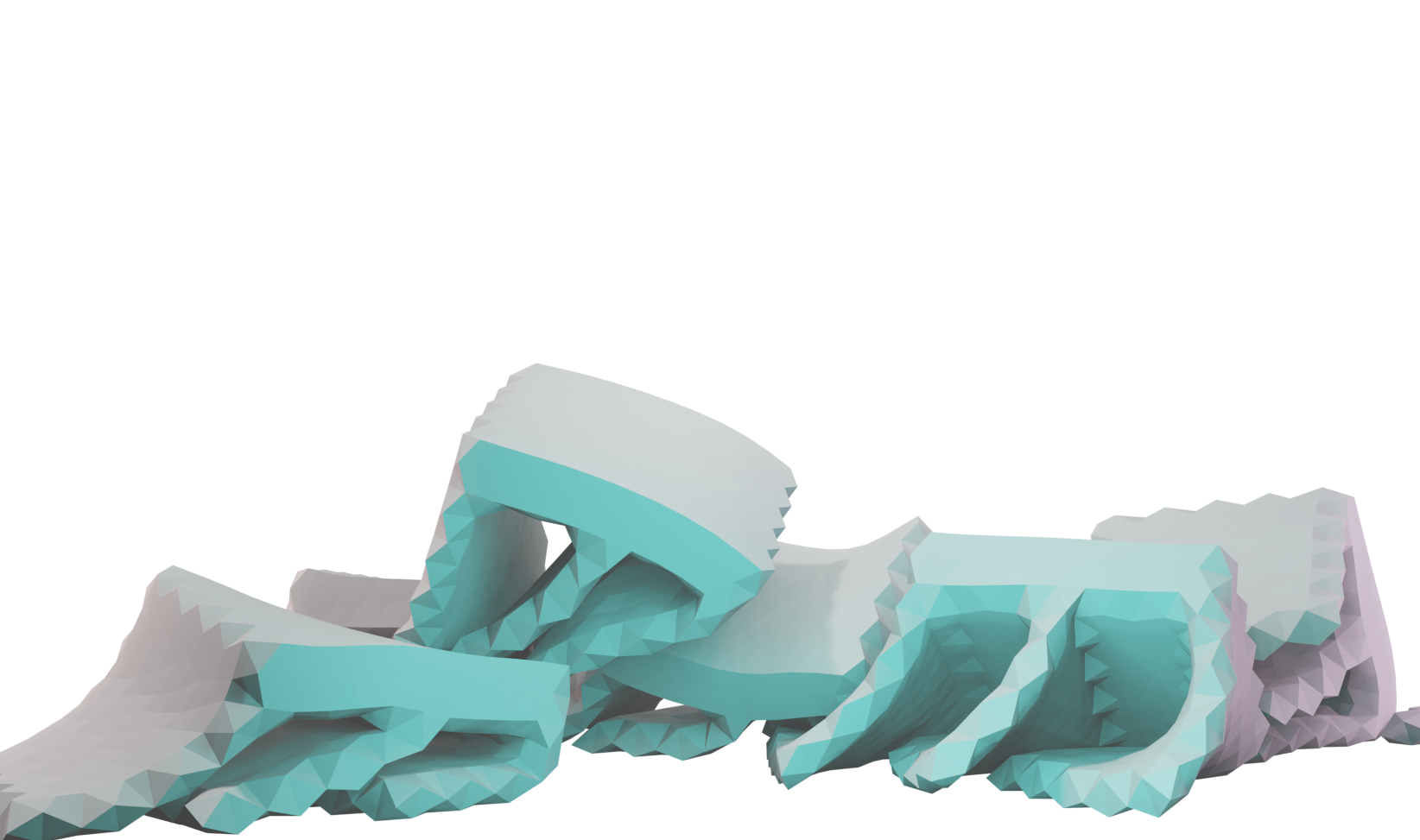}
\end{subfigure}\\ 
\caption{Free fall of 8 "E" models.}
\label{pics:E2}
\end{figure*}

\begin{figure*}[t]
\centering
\includegraphics[width=.245\textwidth]{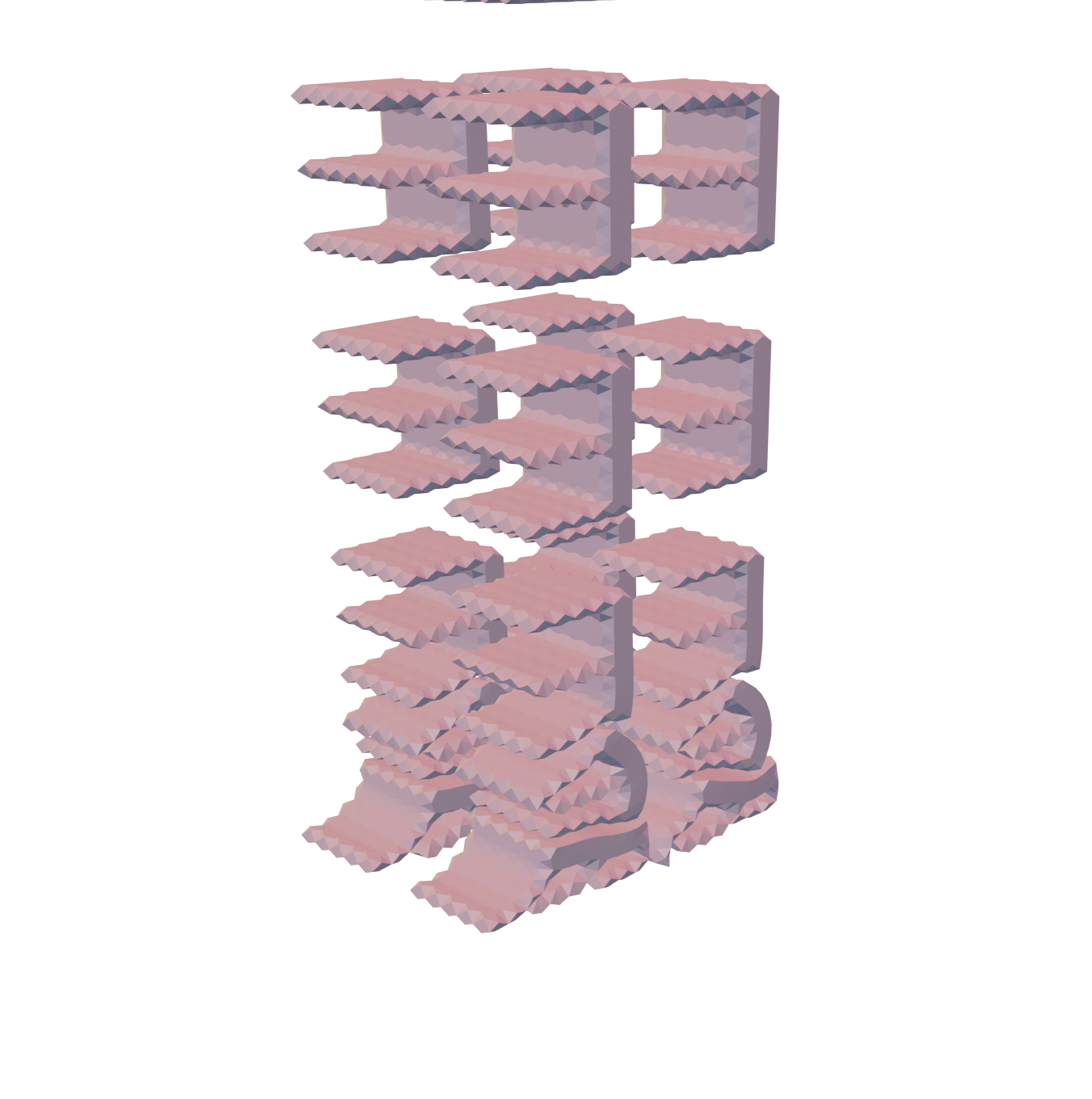}
\includegraphics[width=.245\textwidth]{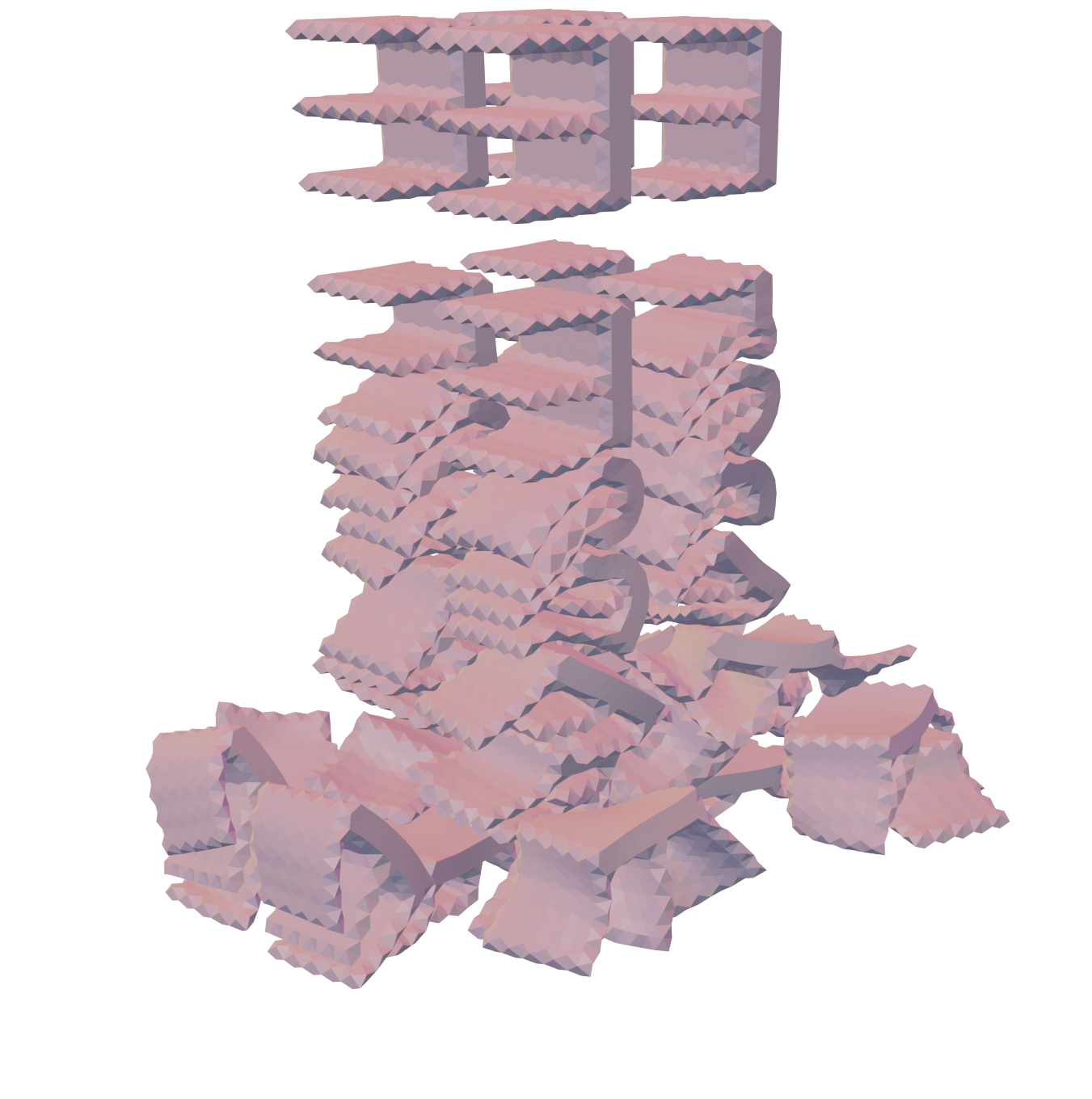}
\includegraphics[width=.245\textwidth]{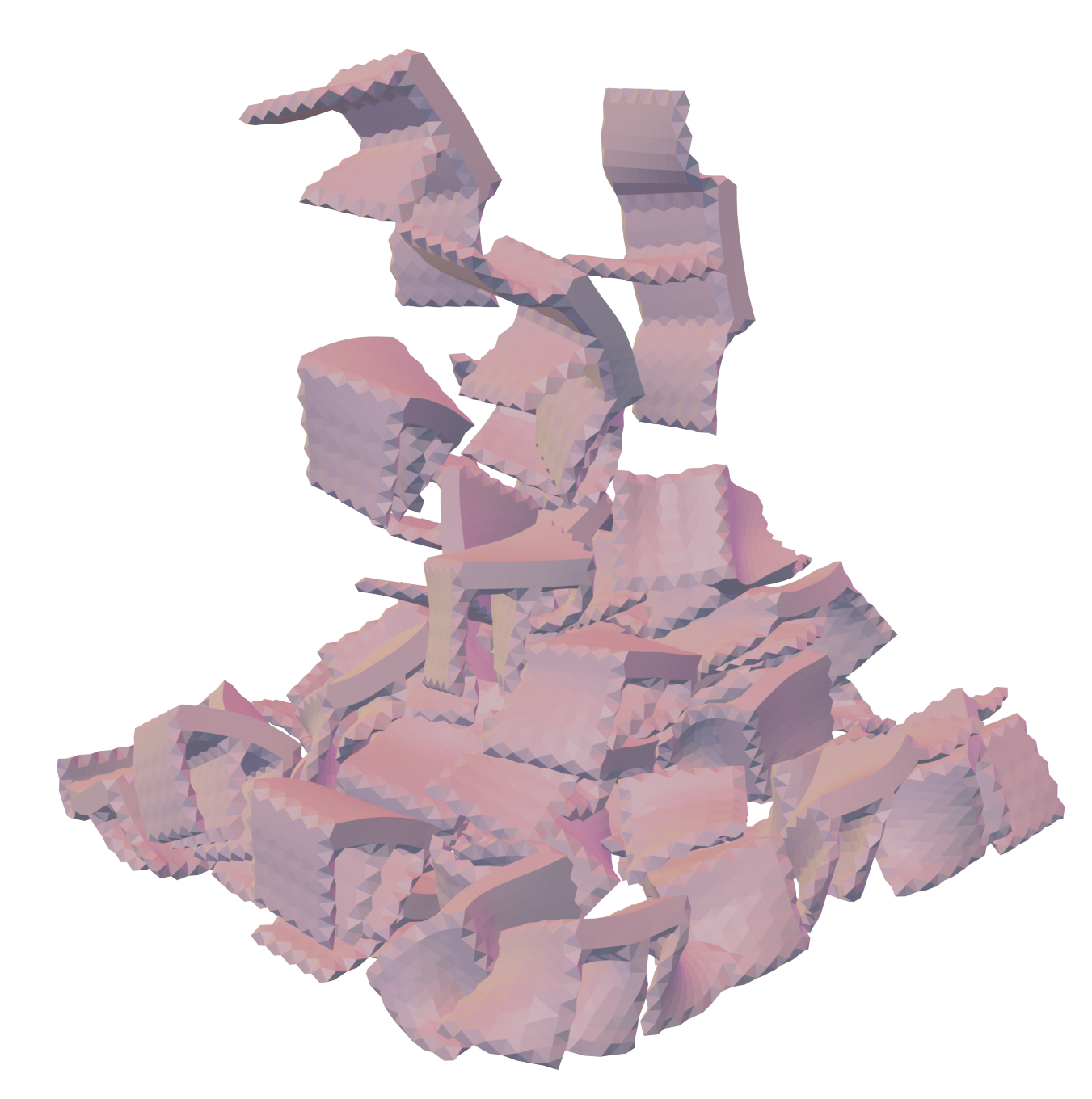}
\includegraphics[width=.245\textwidth]{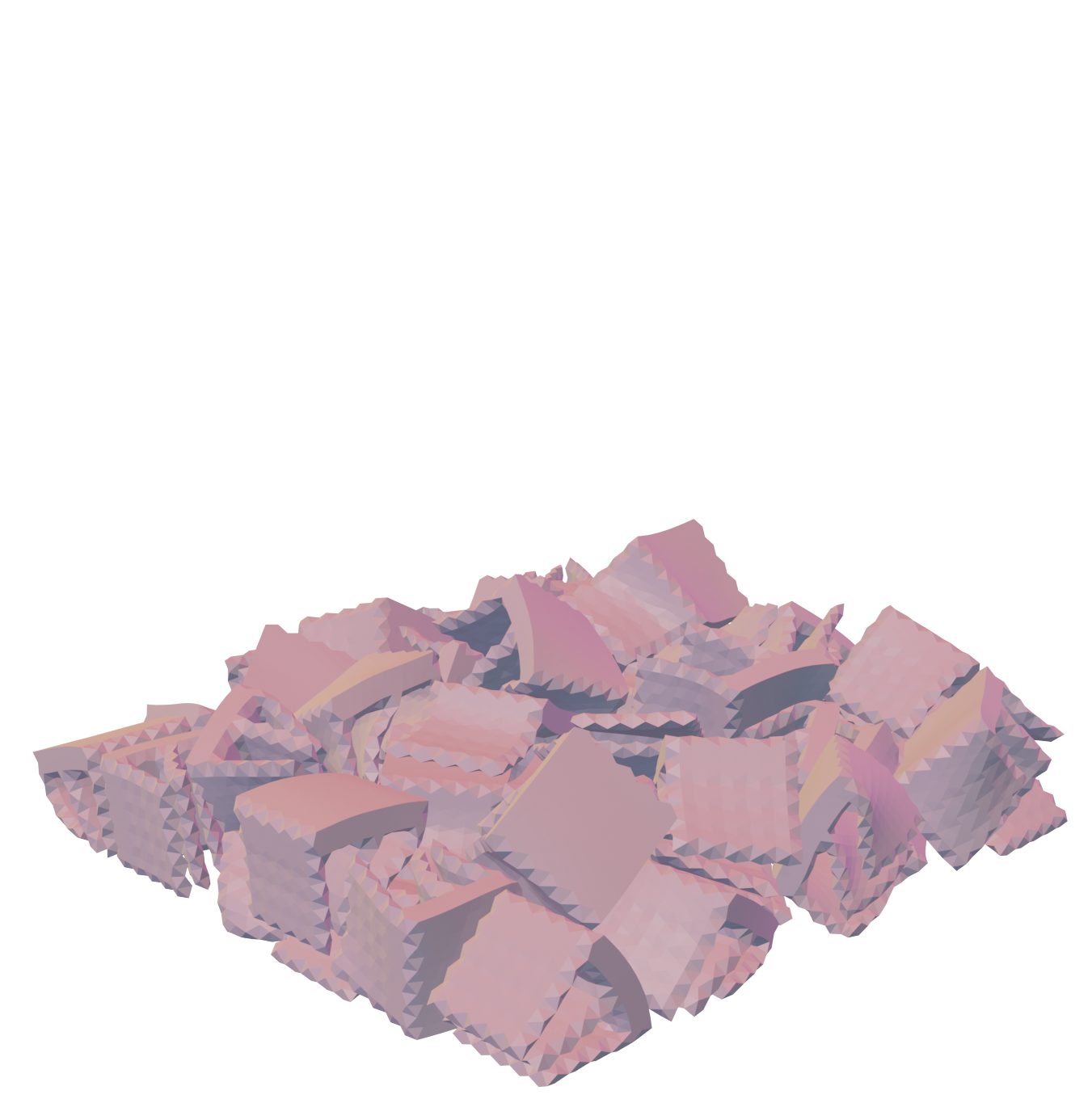}
\caption{Free fall of stacked 48 "E" models.}
\label{pics:E}
\end{figure*}

As illustrated in Table \ref{table}, we evaluate the performance of our method in complex self-collision scenarios. To achieve real-time performance, we employ early stopping of the simulation before it fully converges in the "Squeeze four armadillo," "Four long noodles," and "48 E falling" scenarios. In Figure \ref{pics:armadillo}, we gradually compress four armadillo models from all six directions and subsequently release them, the models exhibit an "explosion" effect upon release, demonstrating correct elastic behavior and the absence of penetrations. In Figure \ref{pics:noodles}, four long strands of noodles simultaneously fall into the bowl with gravitational acceleration 9.8 $m/s^2$. In Figure \ref{pics:E2}, eight instances of the letter "E" undergo free fall toward the ground under a gravitational acceleration of 9.8 $m/s^2$, with a time step of 0.01 s. Figure \ref{pics:E} depicts 48 instances of "E" in free fall, subject to a gravitational acceleration of 0.5 $m/s^2$ and a time step of 0.04 s. As depicted in the corresponding figures, despite the early termination, our method still achieves visually plausible results. Furthermore, we perform discrete collision detection for each frame of all the examples to ensure the absence of penetrations. Futhermore, some unittests \cite{Erleben} are provides in the supplemental video.

\begin{figure*}[htp]
\centering
\includegraphics[width=.245\textwidth]{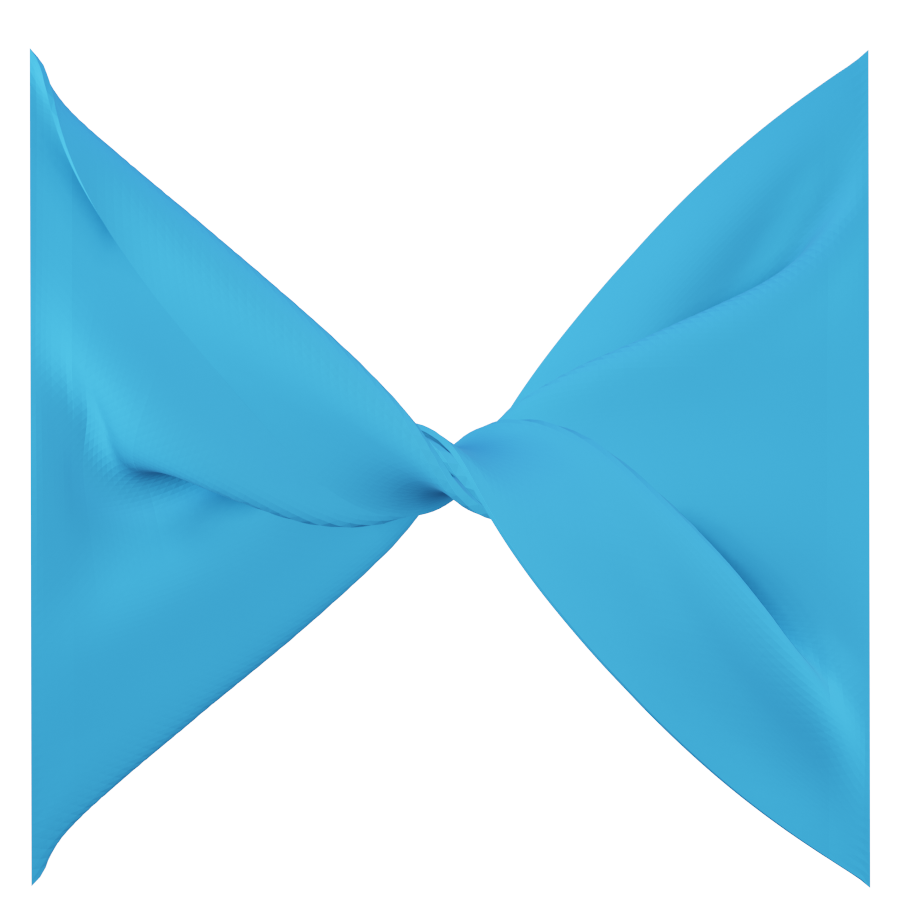}
\includegraphics[width=.245\textwidth]{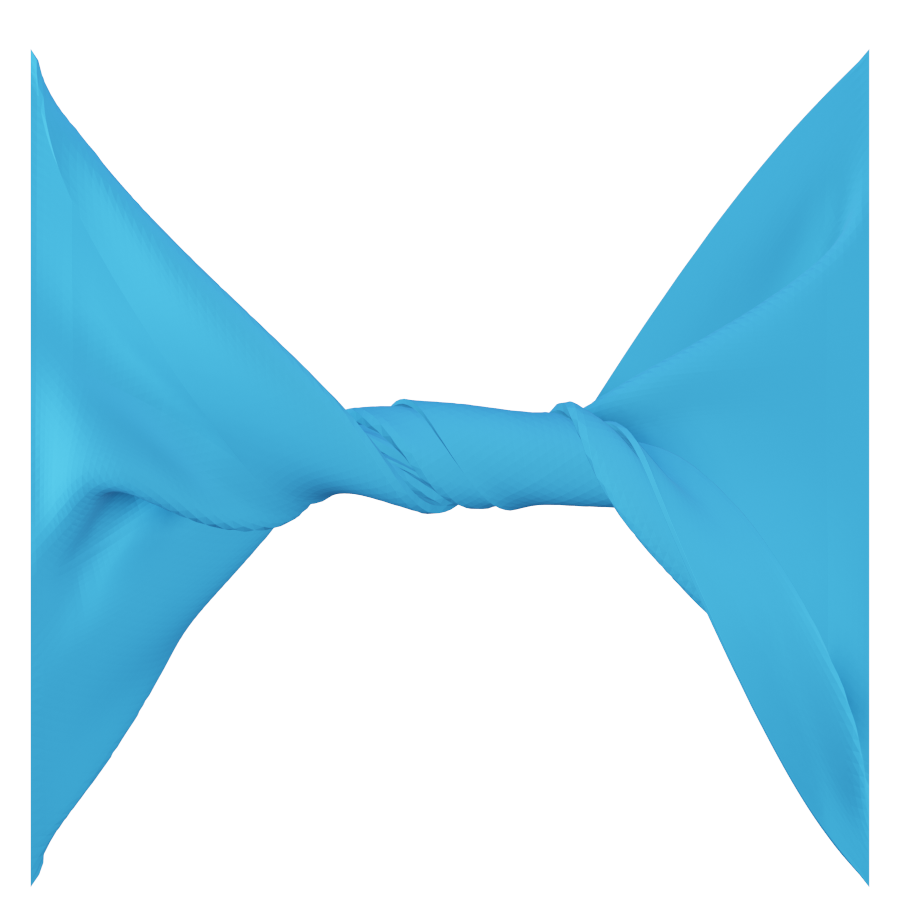}
\includegraphics[width=.245\textwidth]{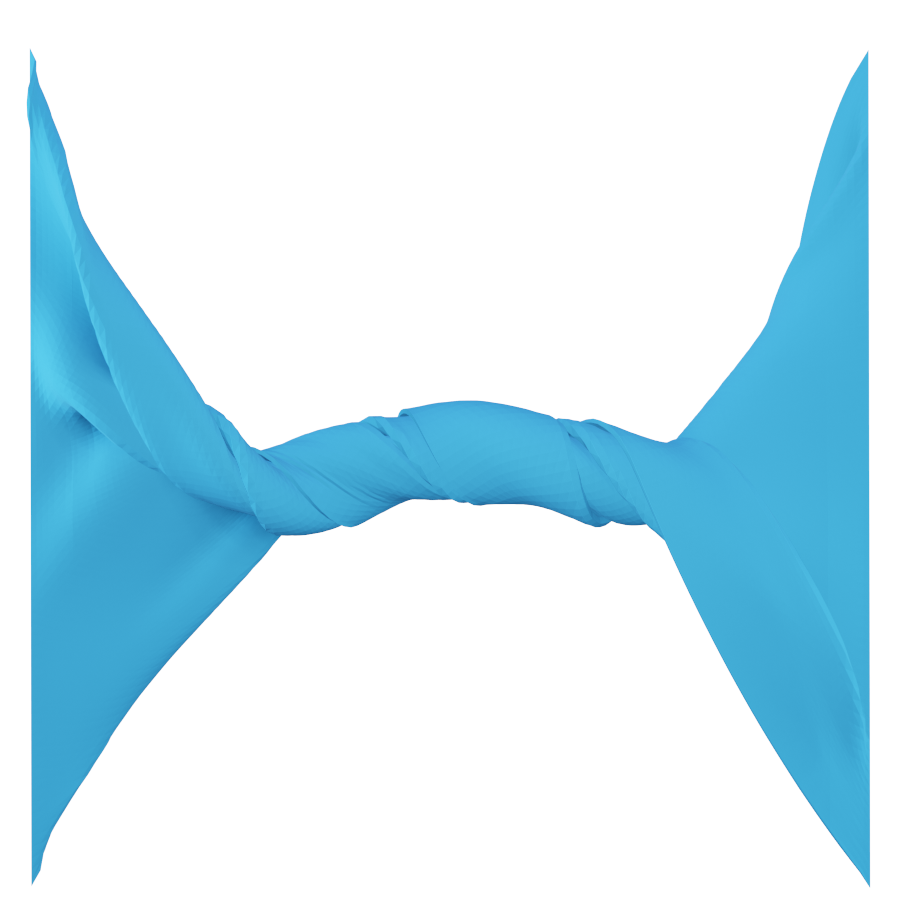}
\includegraphics[width=.245\textwidth]{imgs/mat_blue/250.png}\\
\includegraphics[width=.245\textwidth]{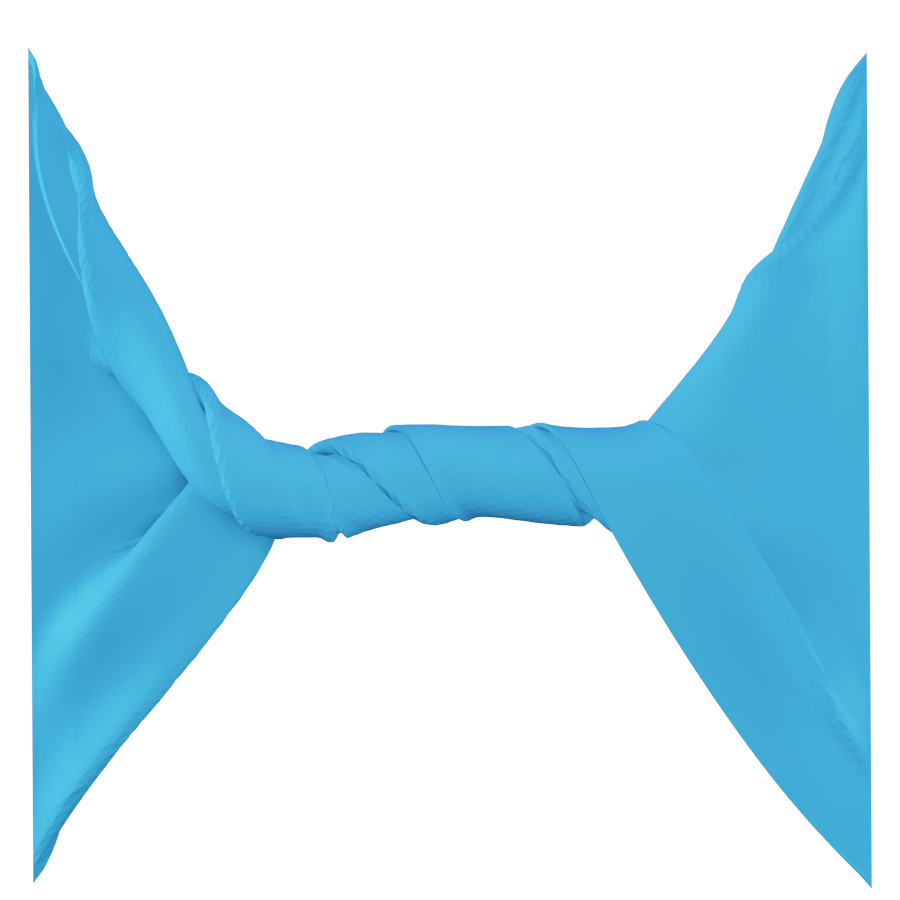}
\includegraphics[width=.245\textwidth]{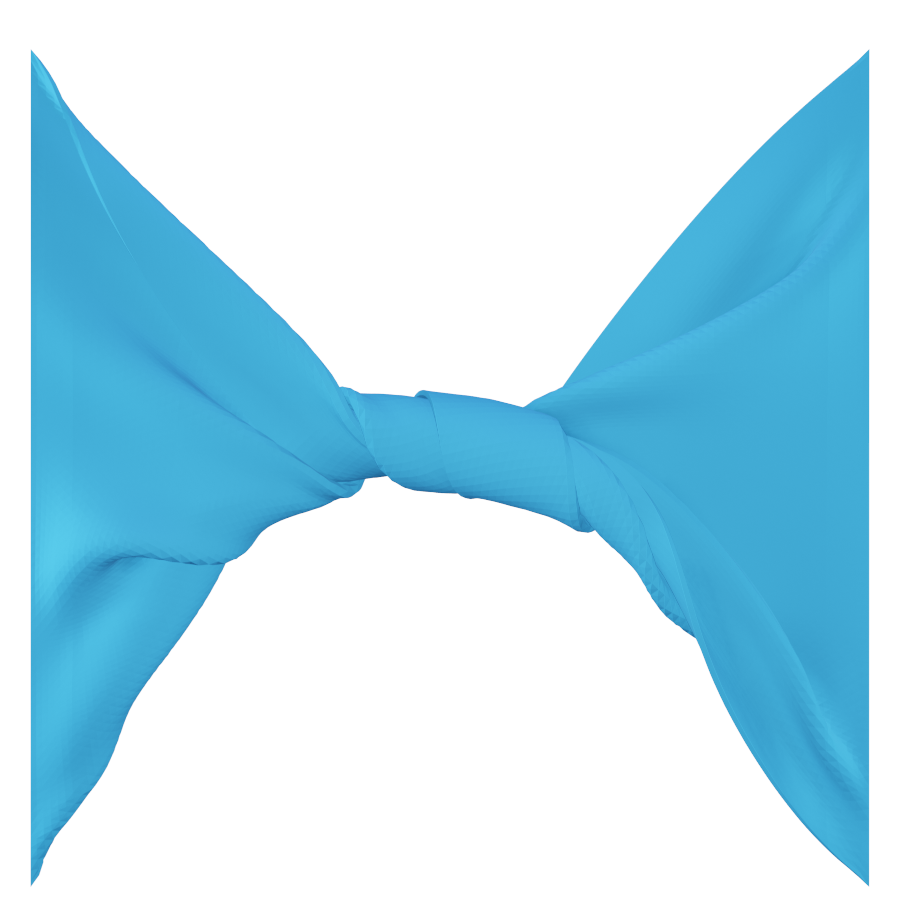}
\includegraphics[width=.245\textwidth]{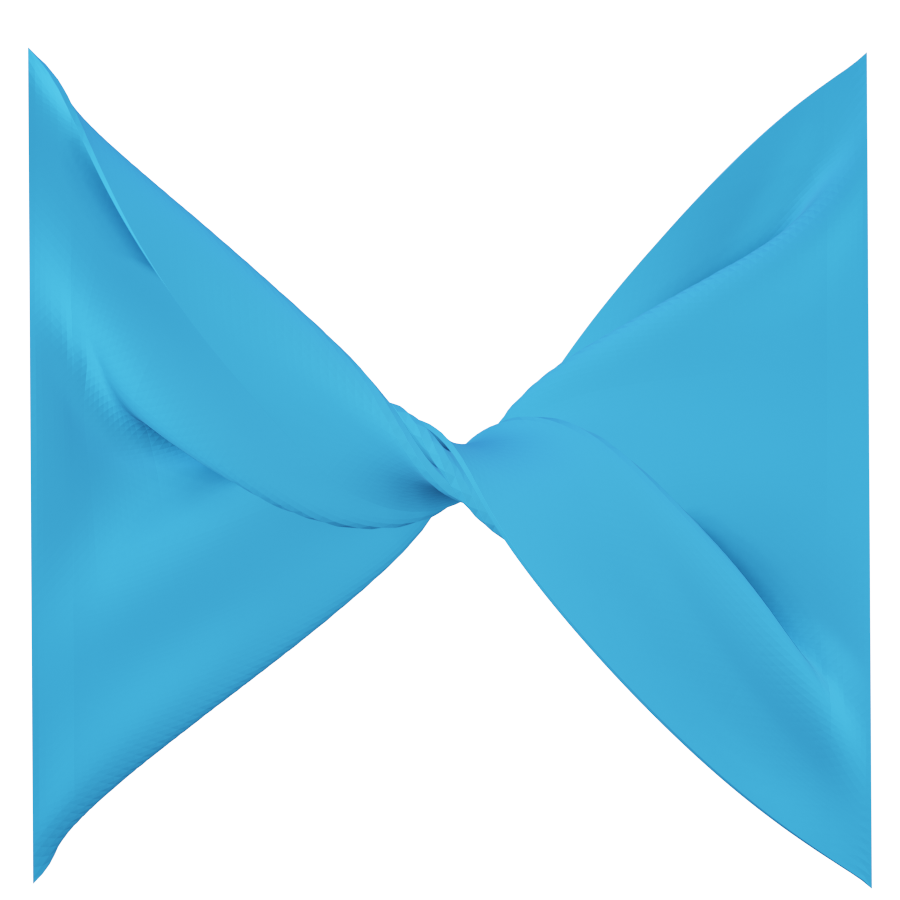}
\includegraphics[width=.245\textwidth]{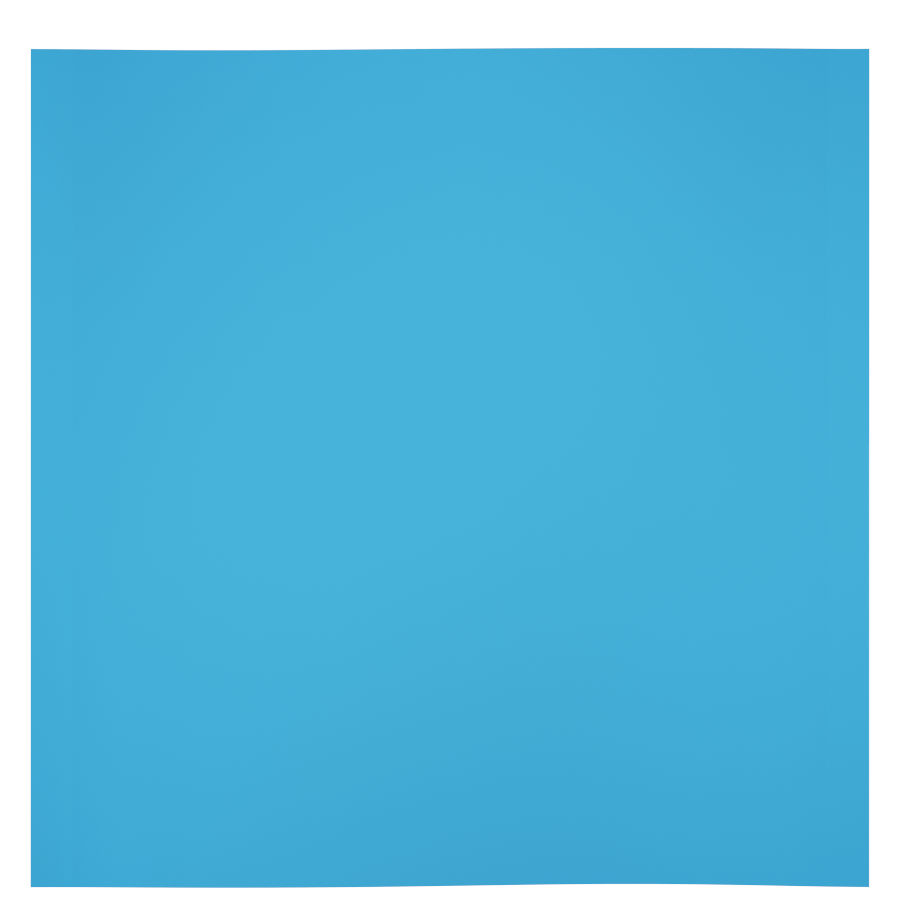}
\caption{Twist mat. In the upper row, the mat is twisted for two rounds, while in the lower row, it is twisted back in the opposite direction. }
\label{pics:twist_mat}
\end{figure*}

In Figures \ref{pics:twist_rods} and \ref{pics:twist_mat}, each side of the rods and mat is twisted with a rotating speed of $72^{\circ}/s$, and notable buckling appears after several rounds of twisting. The mat can still be successfully expanded and recovered by twisting in the opposite direction, demonstrating the absence of penetrations. In the supplementary material, we provide a result of twisting the mat for eight rounds. {As reported in \cite{li2020incremental}, the vanilla CPU IPC requires 33 seconds per time step for the twisting mat demonstrations, whereas our method achieves significantly faster computational speed in comparison}. The state-of-the-art method \cite{li2023second} achieves about 1 FPS speed in similar scenarios. The projective dynamics method \cite{lan2022penetration} displays similar speed, but it does not support general hyperelasticity as our method does.

\section{Conclusion And Limitation}
In this paper, We investigate the application of the nonlinear conjugate gradient method for elastic deformation with barrier contact. Leveraging a fast convergence optimization algorithm, improved line search strategies, simplified Hessian computations, and efficient GPU parallelization, our proposed PNCG method enables real-time simulation of objects comprising over 100K tetrahedra in complex self-collision scenarios. It has the potential for integration with other collision repulsion methods and friction contact.

Since we employ fixed constant values of $\hat{d}$ and $\kappa$ for the IPC potential, and the collision detection module is removed, choosing appropriate values is crucial. While all demonstrations in this paper achieve penetration-free results, in extreme scenarios such as twisting the mat for more rounds, where the mesh undergoes significant deformation, resulting in the dominance of the elastic term, frequently adjusting $\kappa$ is necessary to ensure the updating direction points away from intersections. However, it is challenging to find suitable values that consistently maintain penetration-free results. Similarly, when the object's speed is high and the inertial potential dominates over the elastic and IPC potentials, the step size obtained by Equation (\ref{linesearch}) may cause penetrations. This is the reason why we employ a gravity value of 0.5 in the stacked "48 E falling" example. A robust adaptive barrier stiffness strategy is currently under investigation.


\bibliographystyle{ACM-Reference-Format}
\bibliography{sample-bibliography}

\begin{figure*}[htp]
\centering
\includegraphics[width=.33\textwidth]{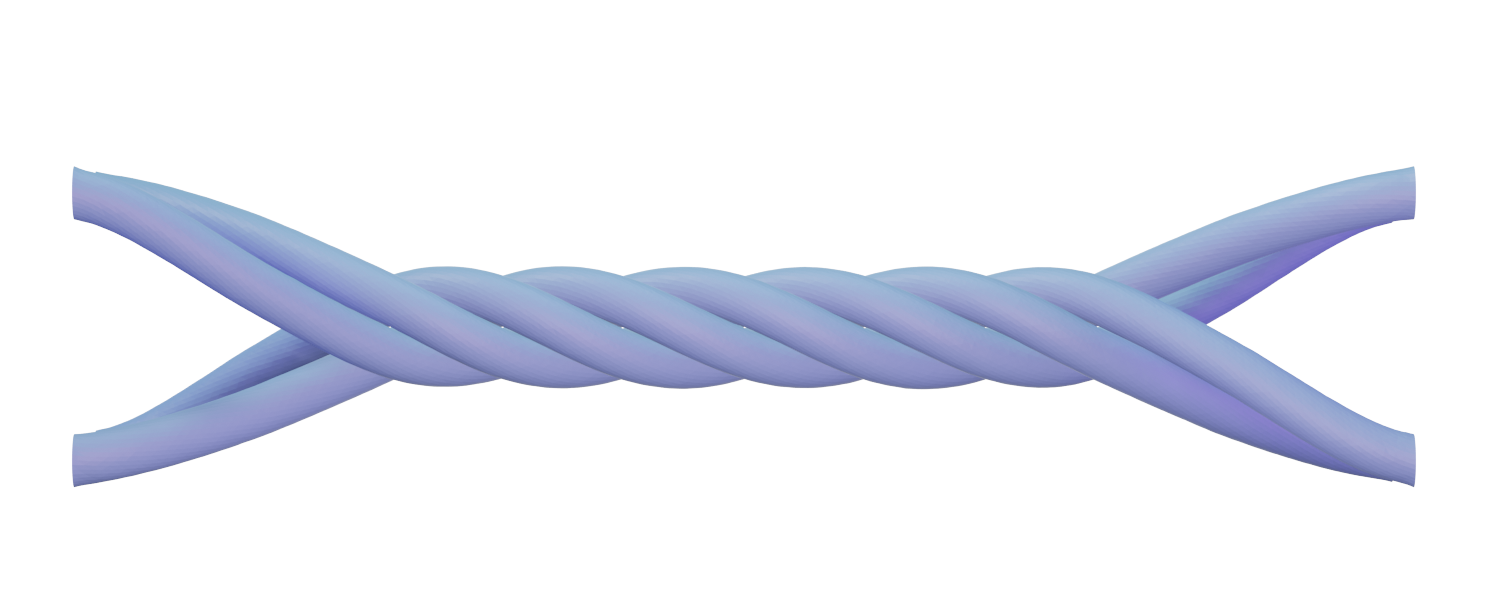}
\includegraphics[width=.33\textwidth]{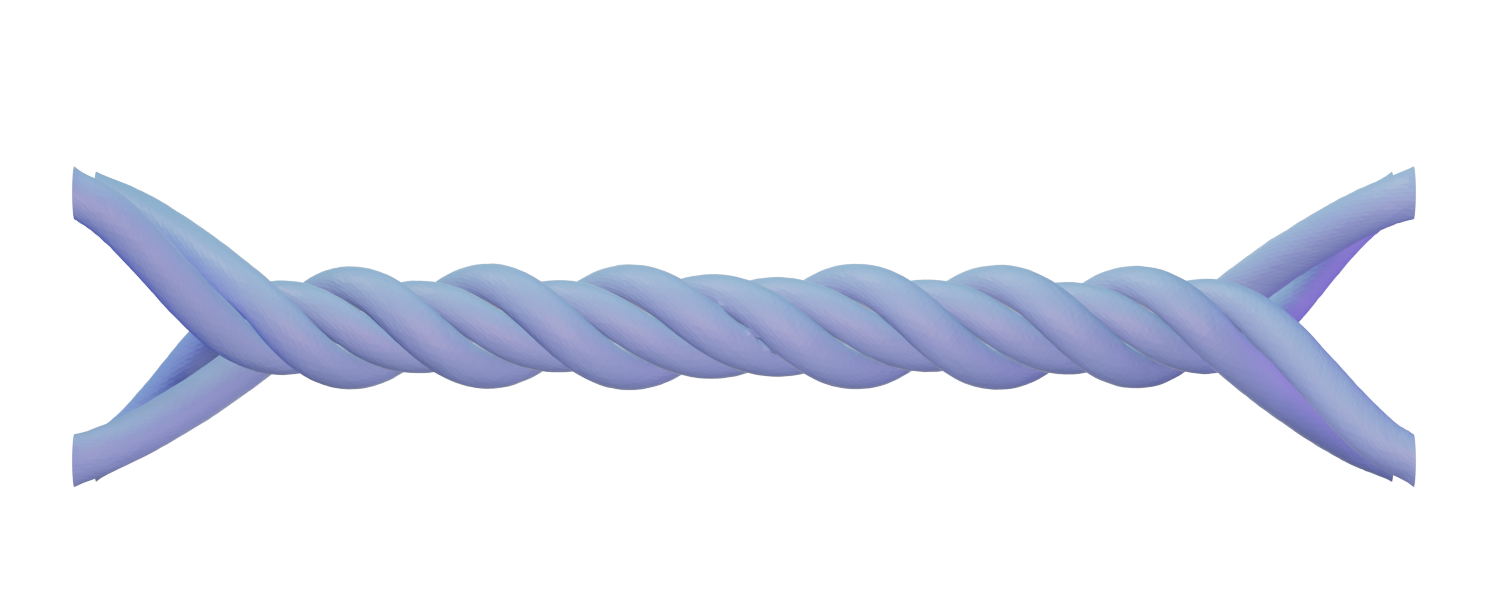}
\includegraphics[width=.33\textwidth]{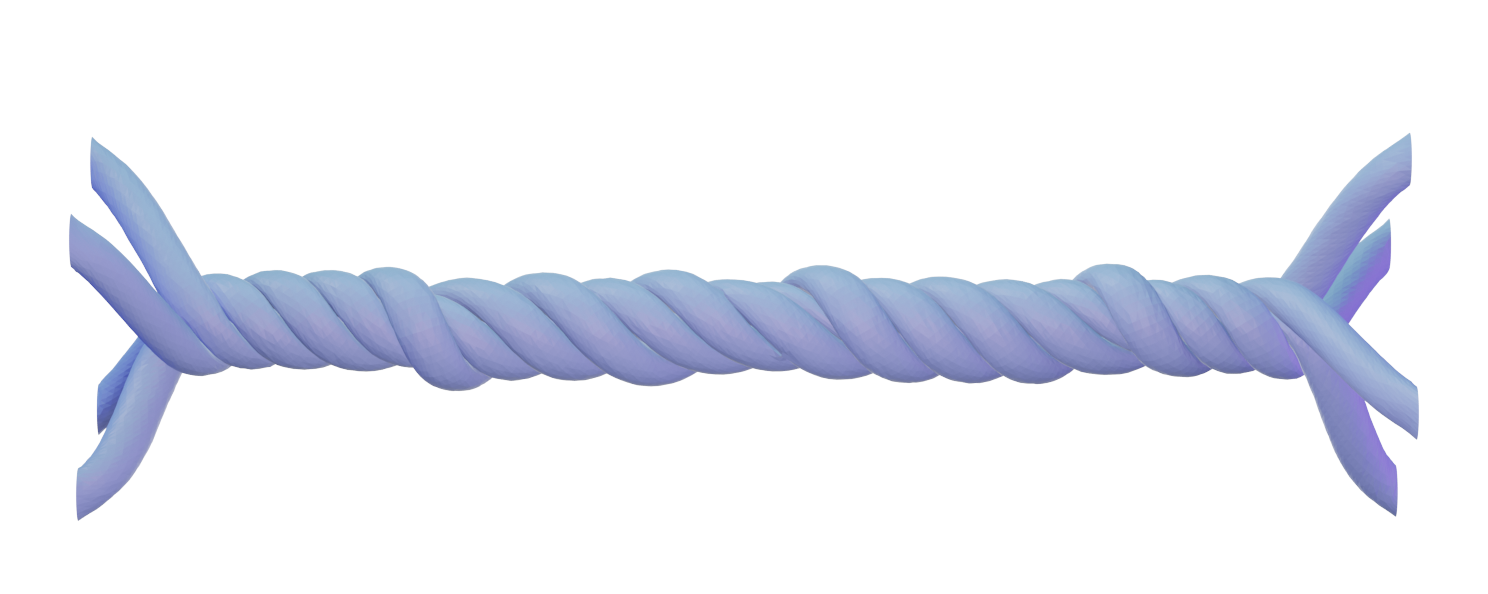}
\caption{Twist rods for 4 rounds.}
\label{pics:twist_rods}
\end{figure*}
\end{document}